\documentclass{cpamart}     
\usepackage{amssymb}
\usepackage{amsmath}

\received{June 2004}       \revised{September 2005}
\yearofpublication{2005}
\volume{000}
\startingpage{1}                      

\authorheadline{Jonathan C. Mattingly \and  \'Etienne Pardoux}
\titleheadline{Malliavin calculus and Navier Stokes equation}

\usepackage{mathptmx}

\newtheorem{theorem}{Theorem}[section]
\newtheorem{lemma}[theorem]{Lemma}
\newtheorem{corollary}[theorem]{Corollary}
\newtheorem{proposition}[theorem]{Proposition}
\theoremstyle{definition}

\theoremstyle{remark}
\newtheorem{remark}[theorem]{Remark}

\newcommand{\hf}        {\frac12}
\newcommand{\ft}        {\frac43}
\newcommand{\N}[1]{\vert\!\vert\!\vert #1 \vert\!\vert\!\vert}
\newcommand{\sol}        {w}
\newcommand{\sgn}        {\textrm{sgn}}
\newcommand{\solf}        {\alpha}
\newcommand{\base}        {e}
\newcommand{\Uf}        {\beta^\phi}
\newcommand{\Idx}        {\mathcal{I}}
\newcommand{\Id}        {I}
\newcommand{\NN}        {\mathbb{N}}
\newcommand{\Ot}        {\Omega_{[0,t]}}
\newcommand{\Zf}        {{\mathcal{Z}_*}}

\newcommand{\R}        {\mathbb{R}}
\newcommand{\DD}        {\mathbb{D}}
\newcommand{\T}        {{\bf T}}
\newcommand{\Z}        {\mathbb{Z}}
\newcommand{\E}        {\mathbb{E}}
\renewcommand{\L}        {\mathbb{L}}
\newcommand{\ccdot}        {\ \cdot \ }
\renewcommand{\P}        {\mathbb{P}}
\renewcommand{\H}        {\mathbb{H}}   

\newcommand{\ip}[2]{\langle#1,#2\rangle}
\newcommand{\eqdef}{\stackrel{\mbox{\tiny def}}{=}}
   
\newcommand{\ONE}{\pmb{1}}
\newcommand{\bJ}{\bar J }
\newcommand{\Mal}{\mathcal{M}}
\newcommand{\hol}{\mathcal{H}}
\newcommand{\OII}{\Omega_{*}'}
\newcommand{\OI}{\hat \Omega_{*}}
\newcommand{\Oc}{\Omega_{c}}

\newcommand{\Os}{\Omega_\sharp'}

\newcommand{\Ospp}{\Omega_\sharp}

\newcommand{\bpf}{\begin{proof}}
\newcommand{\epf}{\end{proof}}
\newcommand{\abs}[1]{|\!\!|\!\!|#1|\!\!|\!\!| }

\begin{document}                        


\title{Malliavin Calculus for the Stochastic  \\2{D} Navier--Stokes Equation}


\author{Jonathan C. Mattingly 
\affil{School of Mathematics, Institute for
  Advanced Studies, Princeton N.J, 08540 USA and Department of
  Mathematics, Duke University, Box 90320, Durham, NC 27708-0320
  USA. Email: jonm@math.duke.edu}
%
%
\\ AND 
\\  \'Etienne Pardoux 
\affil{LATP/CMI, Universit\'e de Provence, 39 rue F. Joliot
  Curie, 13 453 Marseille cedex 13, France. Email: pardoux@cmi.univ-mrs.fr}
}






\begin{abstract}
 We consider the incompressible, two dimensional Navier Stokes
  equation with periodic boundary conditions under the effect of an
  additive, white in time, stochastic forcing. Under mild restrictions
  on the geometry of the scales forced, we show that any finite
  dimensional projection of the solution possesses a smooth, strictly
  positive density with respect to Lebesgue measure.  In particular,
  our conditions are viscosity independent.  We are mainly interested
  in forcing which excites a very small number of modes.  All of the
  results rely on proving the nondegeneracy of the infinite dimensional
  Malliavin matrix.
\end{abstract}

\maketitle   




\section{Introduction}

We consider the movement of a two-dimensional, incompressible fluid
with mean flow zero under periodic boundary conditions. We analyze the
problem using the vorticity formulation of the following form
\begin{equation}
  \label{eq:vorticity}
 \left\{ \begin{aligned}
  &\frac{\partial \sol}{\partial t}(t,x)+  B(\sol,\sol)(t,x) = \nu\Delta
  \sol(t,x) + \frac{\partial W}{\partial t}(t,x) \\
  &\sol(0,x)=\sol_0(x), 
  \end{aligned}\right.
\end{equation}
where $x=(x_1,x_2) \in \T^2$, the two-dimensional torus $[0,2\pi]
\times [0,2\pi]$, $\nu >0$ is
the viscosity constant, $\frac{\partial W}{\partial t}$ is a
white-in-time stochastic forcing to be specified below, and 
\begin{align*}
  B(\sol,\tilde \sol) = \sum_{i=1}^2 u_i(x) \frac{\partial \tilde
  \sol}{\partial x_i}(x)
\end{align*}
where $u=\mathcal{K}(\sol)$. Here $\mathcal{K}$ is the Biot-Savart
integral operator which will be defined next. First we define a
convenient basis in which we will perform all explicit calculations.
Setting $\Z^2_+=\{(j_1,j_2) \in \Z^2 : j_2 > 0\} \cup \{(j_1,j_2) \in
\Z^2 : j_1 > 0, j_2=0\}$, $\Z^2_-=-\Z^2_+$ and $\Z^2_0=\Z^2_+ \cup
\Z^2_-$, we define a real Fourier basis for functions on $\T^2$
with zero spatial mean by
\begin{align*}
  \base_k(x)=
  \begin{cases}
    \sin(k\cdot x) & k \in \Z^2_+\\
    \cos(k\cdot x) & k \in \Z^2_-\ .
  \end{cases}
\end{align*}
We write $\sol(t,x)=\sum_{k \in \Z^2_0} \solf_k(t) \base_k(x)$ for the
expansion of the solution in this basis.  With this notation in the
two-dimensional periodic setting we have the expression
\begin{align}
  \label{eq:BiotSavart}
  \mathcal{K}(\sol)= \sum_{k \in \Z^2_0}
  \frac{k^\perp}{\abs{k}^2}\solf_k \base_{-k},
\end{align}
where $k^\perp=(-k_2,k_1)$ and $\abs{k}^2=k_1^2+k_2^2$.  See for
example \cite{b:MajdaBertozzi02} for more details on the deterministic
vorticity formulation in a periodic domain. We use the vorticity
formulation for simplicity. All of our results can be translated into
statements about the velocity formulation of the problem.

We take the forcing to be of the form
\begin{align}
  \label{eq:W}
  W(t,x)=\sum_{ k \in \mathcal{Z}_*} W_k(t) \base_k(x)
\end{align}
where $\mathcal{Z}_*$ is a finite subset of $\Z^2_0$ and $\{ W_k: k
\in \mathcal{Z}_*\}$ is a collection of mutually independent standard
scalar Brownian Motions on a probability space $(\Omega,\mathcal{F},\P)$. 
The fact
that we force a finite collection of Fourier modes becomes important
starting in Section \ref{sec:hypo}. Up until then the analysis applies
to a force acting on any linearly independent collection of functions
from $\T^2$ into $\R$ which have spatial mean zero. The collection
could even be infinite with a mild summability assumption.

We assume that $\sol_0 \in \L^2=\{ \sol \in L^2(\T^2,\R): \int
\sol\,dx = 0\}$. We will use $\|\ccdot\|$ to denote the norm on $\L^2$
and $\ip{\ccdot}{\ccdot}$ to denote the innerproduct.  We also define
$\H^s=\{ \sol \in H^s(\T^2,\R): \int \sol\,dx = 0\}$. Under these
assumptions, it is standard that $\sol \in C([0,+\infty);{\L^2}) \cap
L^2_{loc}((0,+\infty);{\H}^1)$ \cite{b:Fl94,b:DaZa96,b:MikuleviciusRozovskii02}. 
We
will denote by $\|\cdot \|_s$ the natural norm on $\H^s$ given by
$\|f\|_s=\|\Lambda^s f\|$ where $\Lambda^2=(-\Delta)$.

Our first goal is to prove the following Theorem which will be the
consequence of the more general results given later in the text. In
particular, it follows from Theorem \ref{thm:main}, Theorem
\ref{thm:mainSmooth} and Corollary \ref{c:RussianControl} 
when combined with Proposition
\ref{p:crit}.
\begin{theorem}\label{theo1.1}
  Consider the forcing
  \begin{equation*}
    W(t,x)=  W_1(t)\sin(x_1)+ W_2(t)\cos(x_1) + W_3(t)\sin(x_1+x_2) +
    W_4(t)\cos(x_1+x_2) ,
  \end{equation*}
  then for any $t > 0$ and any finite dimensional subspace $S$ of
  $\L^2$, the law of the orthogonal projection $\Pi \sol(t, \cdot)$ of
  $\sol(t, \cdot)$ onto $S$ is absolutely continuous with respect to
  the Lebesgue measure on $S$. Furthermore, the density is $C^\infty$
  and everywhere strictly positive.
\end{theorem}

A version of Theorem \ref{theo1.1} for Galerkin approximations of
\eqref{eq:vorticity} was one of the main ingredients of the ergodic
and exponential mixing results proven in \cite{b:EMattingly00}. There
the algebraic structure of the nonlinearity was exploited to show that
the associated diffusion was hypoelliptic. Here we use similar
observations on the algebraic structure generated by the vectorfields.
However, new tools are required as there exists little theory applying
Malliavin calculus in an infinite dimensional setting. Relevant
exceptions are \cite{b:HolleyStroock81}, \cite{b:Ocone88}, and
\cite{b:EckmannHairer01b}.

In \cite{b:EckmannHairer01b}, Malliavin calculus was used to
establish the existence of a density when all but a finite number of
degrees of freedom were forced. In contrast to the present paper, the
technique developed there fundamentally required that only a finite
number of directions are unforced. The ideas developed in the present
paper could also likely be applied to the setting of
\cite{b:EckmannHairer01b}.

An essential tool in our approach is a representation of the Malliavin
covariance matrix through the solution of a backward (stochastic)
partial differential equation, which was first invented by Ocone,
see \cite{b:Ocone88}, and which is particularly useful when dealing
with certain classes of SPDEs, since in this case (as opposed to that
of finite dimensional SDEs), the fundamental solution of the
linearized equation cannot be easily inverted. Ocone used that
representation in the case where the original equation is a so--called
``bilinear SPDE'' (that is both the coefficients of ``$dt$'' and
``$dW(t)$'' are linear in the solution). In contrast, we use it in the
case of a nonlinear PDE with additive noise. It seems that these are
the only two cases where Ocone's representation of the Malliavin
matrix through a backward (S)PDE can be used, whithout being exposed
to the trouble of handling a stochastic PDE involving anticipative
stochastic integrals. In Ocone's case, the backward PDE is a
stochastic one, while in our case it is a PDE with random
coefficients. 

There has been a lot of activity in recent years exploring the ergodic
properties of the stochastic Navier-Stokes equations and other
dissipative stochastic partial differential equations. The central new
idea was to make use of the pathwise contractive properties of the
dynamics on the small scales and the mixing/smoothing due to the
stochastic forcing on the larger scales. In \cite{b:Mattingly98b} a
determining modes type theorem (see \cite{b:FoiasProdi67} ) was
developed in the stochastic setting. This showed how controlling the
behavior of a finite number of low modes on a time interval of
infinite length was sufficient to control the entire system.  An
important advance was made concurrently in
\cite{b:BricmontKupiainenLefevere01, b:EMattinglySinai00,
  b:KuksinShirikyan00}, where it was shown that if all of the low
modes were directly forced the system was ergodic.  The first two
covered the case of white in time forcing while the later considered
impulsive forcing.  The assumptions of these papers can be restated as
: the diffusion is elliptic on the unstable subspace of the pathwise
dynamics (see \cite{b:Mattingly03Pre} on this point of view). The
present paper establishes the needed control on the low modes when a
``partial hypoelliptic'' assumption is satisfied.  We show that the
forcing need not excite directly all the unstable modes because the
nonlinearity transmits the randomness to the non--directly excited
unstable directions.  Already, the results of this paper have been
used in an essential way in \cite{b:HairerMattingly04} to prove the
ergodicity of the stochastic Navier Stokes equations under mild,
viscosity independent, assumptions on the geometry of the forcing.

This article is organized as follows: In section
\ref{sec:RepMallMatrix}, we discuss the elements of Malliavin calculus
needed in the paper. In particular, we give an alternative
representation of the quadratic form associated to the Malliavin
matrix. This representation is critical to the rest of the article. In
section \ref{sec:hypo}, we explore the structure of the nonlinearity
as it relates to nondegeneracy of the Malliavin matrix which in turn
implies the existence of a density. In section \ref{sec:nonAdapted},
we prove an abstract lemma on the quadratic variation of non-adapted
processes of a particular form which is the key to the results of the
preceding section. In section \ref{sec:horomander}, we discuss the
relationship to brackets of vector fields and the usual proof of
nondegeneracy of the Malliavin matrix. In doing so we sketch an
alternative proof of the existence of a density. In section
\ref{sec:smooth}, we prove that the density, whose existence is given
in section \ref{sec:hypo}, is in fact $C^\infty$. This requires the
abstract results of section \ref{sec:noSmallEVs} which amount to
quantitative versions of the results in section \ref{sec:nonAdapted}.
Finally in Section \ref{sec:positive}, we prove that the density of the
finite dimensional projections of $\sol(t)$ are everywhere positive
under the same conditions which guarantee smoothness.  We then give a
number of concluding remarks and finish with five appendices
containing technical estimates on the stochastic Navier Stokes
equation. In particular, appendix \ref{sec:HigherMalliavainD} proves
that the solution is smooth in the Malliavin sense and appendix
\ref{sec:LipschitzSupremumEstimates} gives control of the Lipschitz
constants in terms of various quantities associated to the solution.

\section{Representation of the Malliavin covariance matrix}\label{sec:RepMallMatrix}

One way to solve the vorticity equation is by letting
$\sol'(t,x)=\sol(t,x)-W(t,x)$, and solving the resulting PDE with
random coefficients for $\sol'$.  It easily follows from that approach
that for each $t > 0$, there exists a continuous map
\begin{equation*}
  \Phi_t : C([0,t];{\R}^{\mathcal{Z}_*}) \rightarrow {\L^2},
\end{equation*}
such that
\begin{equation*}
  \sol(t)=\Phi_t(W_{[0,t]}).
\end{equation*}
In other words, the solution of equation \eqref{eq:vorticity} can be
constructed pathwise. We shall exploit this in Section
\ref{sec:positive}.

For $k \in \Zf$, $h \in L^2_{loc}(\R_+)$, $t>0$, we define, if it exists, the
Malliavin derivative of $\sol(t)$ in the direction $(k,h)$ as
\begin{align*}
  D^{k,h}\:\sol(t) =\underset{\qquad\quad\varepsilon \ \longrightarrow \ 0}{L^2(\Omega,\L^2)\!-\!\lim}
  \frac{\Phi_t(W+ \varepsilon\:H \base_k)-\Phi_t(W)}{\varepsilon},
\end{align*}
where $H(t)=\displaystyle\int_0^t h(s) ds$. In fact, this convergence
holds pathwise, and it is a Fr\'echet derivative. 

We will show that the above derivative exists, for each $h \in
L^2_{loc}(\R_+)$, and moreover that for each $s \in [0,t]$ and $k \in
\Zf$ there exists a random element $V_{k,s}(t)$ in ${\L^2}$, such that
\begin{equation*}
  D^{k,h}\sol(t)=\int_0^t V_{k,s}(t)h(s) ds.  
\end{equation*}
$V_{k,s}(t)$ is then identified with $D_s^k\sol(t)\eqdef
D^{k,\delta_s}\sol(t)$ and solves equation \ref{2.1} below.
\begin{proposition}\label{prop2.1}
For each $s > 0$ and $k \in \Zf$, the linear parabolic PDE
\begin{equation}\label{2.1}
\left\{
\begin{aligned}
  \frac{\partial}{\partial t} V_{k,s}(t)&= \nu \Delta V_{k,s}(t) -
  B(\sol(t),V_{k,s}(t))- B(V_{k,s}(t), \sol(t)),\:\:t \geq s;\\
  V_{k,s}(s)&=\base_k
\end{aligned}
\right.
\end{equation}
has a unique solution
\begin{equation*}
  V_{k,s} \in C([s,+\infty);{\L^2}) \cap L_{loc}^2([s,+\infty);{\H}^1).
\end{equation*}
\end{proposition}
\bpf
See e.g. Constantin, Foias \cite{b:CoFo88}.
\epf

At times we will consider the linearized equation
\eqref{2.1} with arbitrary initial conditions. We write
$J_{s,t}\phi$ for the solution to \eqref{2.1} at time $t$ with initial
condition $\phi$ at time $s$ less than $t$. In this
notation $V_{k,s}(t)=J_{s,t}\base_k$.

Furthermore, Lemma \ref{l:JJmoments} from the appendix implies
that for  all deterministic
initial conditions $\sol(0)$, $p\geq 1$, $\eta >0$, and $T < \infty$,
\begin{equation*}
\E\sup_{0\leq s\le t\le T}\|V_{k,s}(t)\|^{2p}<
c\exp\Big(\eta\|\sol(0)\|^2 \Big)
\end{equation*}
for some $c=c(\nu,p,T,\eta)$.

Clearly, if $h \in L_{loc}^2 ({\R}_+)$,
\begin{equation*}
  V_{k,h}(t) \eqdef \int_0^t V_{k,s}(t) h(s) ds  
\end{equation*}
is the unique solution in 
$C([0,+\infty);\L^2)\cap L^2_{loc}([0,+\infty);\H^1)$
of the parabolic PDE
\begin{equation}\label{eq:V}
  \left\{
    \begin{aligned}
      \frac{\partial V_{k,h}(t)}{\partial t}&= \nu \Delta V_{k,h}(t) -
      B(\sol(t), V_{k,h}(t))-B(V_{k,h}(t),\sol(t))+ h(t)\base_k,\: t \geq
      0,\\
      V_{k,h}(0)&=0.
    \end{aligned}
  \right.
\end{equation}

It is not hard to see that, in the sense
of convergence in $L^2(\Omega;\L^2)$,
\begin{equation*}
  V_{k,h}(t)=\lim_{\varepsilon \rightarrow 0} \frac{\Phi_t(W+\varepsilon
    H\base_k)-\Phi_t(W)}{\varepsilon} \ .   
\end{equation*}

It then follows that $\sol(t) \in H^1(\Omega, {\L^2})\eqdef \{
X:\Omega \rightarrow \L^2 : \E \|X\|^2, \E\int_0^t\|D^k_sX\|^2ds <
\infty$ for all $k \in \Zf$ and finite $t>0\}$ (see Nualart
\cite{b:Nualart95} page 27 and 62 for more details). Furthermore, its
associated infinite dimensional Malliavin covariance matrix is given
by :
\begin{align}
  \label{eq:defMaliavin}
  \Mal(t) &= \sum_{k \in \Zf} \int_0^t V_{k,s}(t)\otimes
  V_{k,s}(t)ds , \intertext{that is to say it is the operator mapping $\phi \in
    \L^2$ to $ \Mal(t)(\phi) \in \L^2$ given by} \Mal(t)(\phi) &= \sum_{k\in\Zf}
  \int_0^t \ip{V_{k,s}(t)}{\phi} V_{k,s}(t) ds
\end{align}
It follows from Theorem 2.1.2 in Nualart \cite{b:Nualart95} that
Theorem \ref{theo1.1} is a consequence of the fact that for each $\phi
\in {\L^2}$ with $\phi \not=0$,
\begin{equation*}
  \ip{\Mal(t)\phi}{\phi} = \sum_{k \in \Zf} \int_0^t
  \ip{V_{k,s}(t)}{\phi}^2 ds > 0 \mbox{ a.s.} 
\end{equation*}
We now want to give an alternative representation of this quantity,
using a backward PDE which is the adjoint of equation (\ref{eq:V}).
\begin{proposition}\label{prop2.2}
For each $t > 0$, $\phi \in {\L^2}$, the linear backward parabolic PDE
\begin{equation}\label{2.2}
\left\{
\begin{aligned}
  \frac{\partial}{\partial s} U^{t,\phi}(s)+\nu\:\Delta U^{t,\phi}(s)+
  B(\sol(s),U^{t,\phi}(s))-C(U^{t,\phi}(s),\sol(s))=0,\\
  \hfill  0 \leq s \leq t;\\
  U^{t,\phi}(t)=\phi,
\end{aligned}
\right.
\end{equation} 
has a unique solution
\begin{equation*}
  U^{t,\phi} \in C ([0,t];{\L^2}) \cap L^2([0,t) ; {\H}^1).
\end{equation*}
Here $C\big(\ccdot,\sol(s)\big)$ is the $\L^2$-adjoint of the time-dependent,
linear operator $B(\ccdot,\sol(s))$ and thus is defined by the
relation $\ip{B(u,\sol(s))}{v}=
\ip{C(v,\sol(s))}{u}$.
\end{proposition}
\bpf Same argument as in Proposition \ref{prop2.1} \epf 

As before Lemma \ref{l:JJmoments} from the appendix implies that
there exist a positive constant $\eta$ so that for all deterministic
initial conditions $\sol(0),\phi \in \L^2$, $p\geq 1$ and $T < \infty$
\begin{equation*}
\E\sup_{0\leq s\le t\le T}\|U^{t,\phi}(s)\|^{2p}<
c\|\phi\|^{2p}\exp\Big(\eta\|\sol(0)\|^2 \Big). 
\end{equation*}
for some $c=c(\nu,p,T,\eta)$.

\begin{proposition}\label{prop2.3}
For each $k \in \Zf$ and $\phi \in \L^2$, the function
\[ 
r \to \ip{V_{k,s}(r)}{U^{t,\phi}(r)}
\]
from $[s,t]$ into ${\R}$ is constant.
\end{proposition}
\bpf We first show that this mapping belongs to
$W^{1,1}(s,t;\R)$. It is clearly a continuous function and moreover
\begin{equation*}
  V, U \in C\big([s,t];{\L^2}\big) \cap L^2\big((s,t);{\H}^1\big),
\end{equation*}
hence by interpolation
\begin{equation*}
  V, U \in L^4\big((s,t),{\H}^{\hf}\big) \cap L^2\big((s,t);{\H}^1\big).  
\end{equation*}
Since we also know that 
\begin{equation*}
  \sol \in L^4((s,t);{\H}^{\hf}) \cap L^2((s,t);{\H}^1),  
\end{equation*}
we know that the products $\|U\|_{\hf} \|\sol\|_1,\|V\|_{\hf}
\|\sol\|_1$, $\|w\|_{\hf} \|U\|_1$ and $\|w\|_{\hf} \|V\|_1$ all belong
to $L^{\ft}\big((s,t);\R\big)$.  On the other hand, it follows from the fact that
$\|\mathcal{K}U\|_1=\|U\|$ and  estimate
(6.10) in Constantin, Foias \cite{b:CoFo88} that
  \begin{xalignat*}{2}
    |\ip{C(U,\sol)}{\psi}|&\leq c\|U\|_{\hf} \|\sol\|_1 \|\psi\|_{\hf} &
    |\ip{B(V,\sol)}{\psi}|&\leq c\|V\|_{\hf} \|\sol\|_1  |\|\psi\|_{\hf}\\
    |\ip{B(\sol,U)}{\psi}|&\leq c\|\sol\|_{\hf} \|U\|_1 \|\psi\|_{\hf}&
    |\ip{B(\sol,V)}{\psi}|&\leq c\|\sol\|_{\hf} \|V\|_1 \|\psi\|_{\hf} .
  \end{xalignat*}
  Hence we conclude that $C(U,\sol)$, $B(V,\sol)$, $B(\sol,U)$ and
  $B(\sol,V)$ all belong to $L^{\ft}((s,t);{\H}^{-\hf})$. From
  \eqref{2.1} and \eqref{2.2}, we that that both 
  $\frac{d\ }{dt}V$ and $\frac{d\ }{dt}U$
  consist of three terms. The first belongs to $L^2((s,t);{\H}^{-1})$
  and the last two to $L^{\ft}((s,t);{\H}^{-\hf})$. Hence
  $\ip{V_{k,s}(r)}{\frac{d\ }{dr}U^{t,\phi}(r)}$ and
  $\ip{\frac{d\ }{dr}V_{k,s}(r)}{U^{t,\phi}(r)}$ are in
  $L^1(s,t; \R)$ and the statement that
\begin{equation*}
  r \to \ip{V_{k,s}(r)}{U^{t,\phi}(r)}  
\end{equation*}
is a.e. differentiable then follows a variant of Theorem 2, Chapter
18, section 1 in Dautray, Lions \cite{b:DautrayLions88}. Moreover, for
almost every $r$
\begin{equation*}
  \begin{split}
    \frac{d}{dr}\ip{V(r)}{U(r)}&=
    \ip{A(\sol(r))V(r)}{U(r)}-\ip{V(r)}{A^*(\sol(r)) U(r)}\\ 
    &=0
  \end{split}
\end{equation*}
where $A(\sol(t))$ is the linear operator on the right handside of
\eqref{2.1} and $A^*(\sol(t))$ is its $\L^2$-adjoint. The result
follows.  \epf

We can now rewrite the Malliavin covariance matrix using $U$ in place
of $V$. For a fixed $\phi$, this is an improvement as $U^{t,\phi}(r)$
is a single solution to a PDE while $V_{k,t}(r)$ is a continuum of
solutions indexed by the parameter $s$.
\begin{corollary}\label{cor2.4}
  For any $\phi \in {\L^2}$,
  \begin{equation*}
    \ip{\Mal(t) \phi}{\phi} = \sum_{k\in \Zf} \int_0^t \ip{\base_k}{U^{t,\phi}(s)}^2 ds.  
  \end{equation*}
\end{corollary}
\bpf This follows from the fact that for any $0 \leq s \leq t$ and $k
\in \Zf$,
\[\ip{V_{k,s}(s)}{U^{t,\phi}(s)} = \ip{V_{k,s}(t)}{U^{t,\phi}(t)}\]
i.e.
\begin{equation*}
  \ip{\base_k}{U^{t,\phi}(s)}=\ip{V_{k,s}(t)}{\phi}
\end{equation*}
\epf

From this corollary, one immediately deduces the following result.
\begin{corollary}\label{l:baseCase} 
   Assume that for some fixed $\phi \in \L^2$,
  \begin{equation*}
    \ip{\Mal(t)\phi}{\phi} =0
  \end{equation*}
  on a subset $\Omega_1$ of $\Omega$. Then for
  all $k \in \Zf$ and $s\leq t$, $\ip{\base_k}{U^{t,\phi}(s)}=0$ on
  $\Omega_1$. In particular, $\ip{\base_k}{\phi}=0$.
\end{corollary}

\section{Hypoellipticity}
\label{sec:hypo}

\subsection{Final Assumptions and  Main Existence Result}

We define $\mathcal{Z}_0$ to be the symmetric part of the forcing
set $\Zf$ given by $\mathcal{Z}_0= \Zf \cap (-\Zf)$ and then the
collection
\begin{multline*}
  \mathcal{Z}_n= \big\{\ell+j\in \Z^2_0 :
  j \in \mathcal{Z}_0,\; \ell \in  \mathcal{Z}_{n-1}
 \text{with } \ell^\perp\cdot j \not =0,\; \abs{j}\not =\abs{\ell} \big\}
\end{multline*}
and lastly, 
\begin{equation*}
  \mathcal{Z}_\infty=\bigcup_{n=1}^\infty \mathcal{Z}_n .
\end{equation*}
Notice that the above union starts at one and that the $\mathcal{Z}_n$ are
symmetric in that $\mathcal{Z}_n=-\mathcal{Z}_n$. This follows by recurrence,
starting with
$\mathcal{Z}_0=-\mathcal{Z}_0$. We are mainly
concerned with the case where $\mathcal{Z}_0=\Zf$ as this corresponds
to noise which is stationary in $x$. We now can state the main
theorem. Defining
 \begin{equation}
   \label{eq:S}
  S_0=\mathrm{Span}\biggl(\base_k: k \in\Zf\biggr);\; S_n=\mathrm{Span}\biggl( \base_k: k \in \bigcup_{j=1}^n\mathcal{Z}_j
     \cup \Zf\biggr),\; n \in \{1,2,\dots,\infty\},
\end{equation}
we have the following result which implies the first part of Theorem \ref{theo1.1} in
the case when $S_\infty=\L^2$.
\begin{theorem}\label{thm:main}
  For any $t > 0$ and any finite dimensional subspace $S$ of $S_\infty$,
  the law of the orthogonal projection $\Pi \sol(t, \cdot)$ of
  $\sol(t, \cdot)$ onto $S$ is absolutely continuous with respect to
  the Lebesgue measure on $S$.
\end{theorem}

The above result guarantees an absolutely continuous density on finite
dimensional subsets of $S_\infty$. However, it does not imply the lack
of density for other subsets as in constructing $S_\infty$ we have
only used part of the available information. In the proof below, it
will become clear that we make use only of the directions generated by
frequencies where both the $\sin$ and $\cos$ are stochastically
forced. We do this in the name of simplicity and utility. Verifying
any more complicated condition was
difficult. However, as translation invariance implies that both the
$\sin$ and $\cos$ mode of a given frequency are forced, it seems a
reasonable compromise.  In the end, we are primarily interested in
producing conditions which give insight as to how the nonlinearity
spreads the randomness. In particular we now give an easy Proposition
which in conjunction with Theorem \ref{thm:main} proves the first part of
Theorem
\ref{theo1.1} given in the introduction. After that we will quote a
more general result from \cite{b:HairerMattingly04} which is proven
using similar ideas to those below.
\begin{proposition}\label{p:crit}
  If  $\{ (0,1),(1,1)\} \subset \mathcal{Z}_0$ then  $S_\infty=\L^2$.
\end{proposition}
\bpf Clearly by adding and subtracting the vectors $(0,1)$ and
$(1,1)$, one can generate all of $\Z_0^2$. The only question is whether
the conditions $\ell^\perp \cdot j\neq 0$ and $\abs{j}\neq\abs{\ell}$ from
the definition of $\mathcal{Z}_n$ ever create a situation which blocks
continuing generating the lattice. First notice that
$(1,0)=(1,1)-(0,1)$, $(-1,0)=(0,1)-(1,1)$ and $(-1,1)=(-1,0)-(1,1)+(0,1)+(1,1)$ and if
these moves are made from left to right none of the restrictions are
violated. Hence all of the vectors of the same length as the vectors
in $\mathcal{Z}_0$ can be reached and henceforth the restriction
$\abs{j}\neq\abs{\ell}$ will not be binding.  The requirement that $\ell^\perp
\cdot j\neq 0$ does not cause a problem. The line $\ell^\perp \cdot
(0,1) =0$ can be approached from above and does not obstruct
generating the rest of the lattice. The line $\ell^\perp \cdot (1,1)
=0$ does separate the lattice. However we know we can reach $(1,0)$
and have the point $(0,1)$ to start with. Hence we can reach all of
the points on either side of the line. Those on the line can be reached
from points on either side. \epf

It is clear from the preceding lemma that many other choices of
forcing will also lead to $S_\infty=\L^2$. For instance if $\{
(1,0),(1,1)\} \subset \mathcal{Z}_0$ then $S_\infty=\L^2$. It is also
interesting to force a collection of modes distant from the origin and
allow the noise to propagate both up to the large scales and down to
still smaller scales. We now give a simple proposition giving sufficient
conditions in such a setting. 
\begin{proposition}
  Let $M,K \in \NN$ with $M,K > 2$ and $\abs{M-K}> 2$. Then if $\{
  (M+1,0),(M,0),(0,K+1),(0,K) \} \subset \mathcal{Z}_0$, 
  $S_\infty=\L^2$.
\end{proposition}
\bpf The idea is to use $(M+1,0)-(M,0)=(1,0)$,$(M,0)-(M+1,0)=(-1,0)$,
$(0,K+1)-(0,K)=(0,1)$ and $(0,K)-(0,K+1)=(0,-1)$ in order to generate
the whole lattice. The only difficulty could be the above restrictions. 
The restrictions of the form $\ell^\perp
\cdot j\neq0$ only prevent applying $(M+1,0)-(M,0)$ or $(M,0)-(M+1,0)$
to points on the $x$-axis and $(0,K)-(0,K+1)$ and $(0,K+1)-(0,K)$ to
points on the $y$-axis. However this is not a serious restriction as
all of the points on the $x$-axis can be reached by moving down from
above and all of the points on the $y$-axis can be reached by moving
horizontally. Furthermore the $y$-axis can be crossed by using
strictly horizontal moves. The only remaining restriction is the
points $k \in \Z_0^2$ with $\abs{k}\in \{K,K+1,M,M+1\}$. For example
assume that $\abs{k}=K$ and one wanted to move to the left by applying
$(0,K)-(0,K+1)$. While this direct move is illegal, one can accomplish
the same effect by moving up then left and finally down. The
requirement that $\abs{M-K} >2$ ensures that we will not be blocked from
moving up using $(M+1,0)-(M,0)$ given that $k=\abs{K}$. Once we have moved
up, we will be free to move to the left and then back down. The other
cases are analogous. \epf

Guided by these results, in \cite{b:HairerMattingly04} the following is
proven:
\begin{proposition} One has  $S_\infty=\L^2$ if and only if:
  \begin{enumerate}
  \item Integer linear combinations of elements of $\mathcal{Z}_0$
  generate $\Z^2_0$.
\item There exist at least two elements in  $\mathcal{Z}_0$  with
  unequal euclidean norm.
  \end{enumerate}
\end{proposition}
This gives a very satisfactory characterization of the setting when
$S_\infty=\L^2$ which is the case of primary interest.

\subsection{Proof of Theorem \ref{thm:main}}
       
Since we already know that $\sol(t) \in H^1(\Omega,{\L^2})$, Theorem
\ref{thm:main} follows from Theorem 2.1.2 in Nualart \cite{b:Nualart95} and the
fact that for any $\phi \in S_\infty$, $\phi \not=0$,
\begin{equation*}
\ip{\Mal(t) \phi}{\phi}\,>0\:\mbox{  a.s.}  
\end{equation*}
Hence to prove Theorem \ref{thm:main}, it suffices to show that
\begin{proposition}\label{lem3.1}
There exists a subset $\Omega_1 \subset \Omega$ of full measure so
that on $\Omega_1$ if $\ip{\Mal(t) \phi}{\phi} = 0$ for some $\phi \in
\L^2$ then $\Pi_\infty\phi=0$ where $\Pi_\infty$ is the
$\L^2$-orthogonal projection onto $S_\infty$.
\end{proposition}
Notice that Proposition \ref{lem3.1} is equivalent to 
\begin{equation}\label{eq:noDeg}
       \P\Big( \bigcap_{\substack{\phi \in \L^2\\ \Pi_\infty \phi \neq
      0}} \big\{ \ip{\Mal(t)\phi}{\phi} > 0 \big\}  \Big)=1 \ .
\end{equation}
To prove the proposition we need to better understand the structure
of the equations. To this end, we now write the equations for the
spatial Fourier coefficients of $U$ and $\sol$ to better expose the
interactions between the systems various degrees of freedom. In this and the
general structures of the nonlinearity exploited, we follow E,
Mattingly \cite{b:EMattingly00}; however, the tact of the analysis is
different. (We also take the chance to correct a small error in
\cite{b:EMattingly00}. There the summation was restricted to modes in
the first quadrant, when it should have ranged over the entire
upper-half plane.) Again setting $\sol(t,x)=\sum_{k\in \Z_0^2}
\solf_k(t)\base_k$, we have for $\ell \in \Z^2_0$
\begin{multline}
  \label{eq:fourierAlpha}
  \frac{d}{dt}\solf_\ell(t)+\nu\abs{\ell}^2\solf_\ell(t) +
  \frac12\sum_{(j,k,\ell)\in \Idx_+} c(j,k)\solf_j(t)
  \solf_k(t)\\-\frac12 \sum_{(j,k,\ell)\in \Idx_-}
  c(j,k)\solf_j(t) \solf_k(t) =\ONE_{\Zf}(\ell)  \frac{d}{dt}W_\ell(t)
\end{multline}
where $\ONE_{\Zf}$ is the indicator function of $\Zf$, $c(j,k)=\frac12
(j^\perp\cdot k) (\abs{j}^{-2}- \abs{k}^{-2})$ and 
\begin{align*}
    \Idx_+=& \big\{  (j,k,\ell) \in  (\Z^2_+,\Z^2_-,\Z^2_+)\cup (\Z^2_-,\Z^2_+,\Z^2_+) \cup(\Z^2_+,\Z^2_+,\Z^2_-)\;\big|\; k+j+\ell=0 \big\} \\&
  \cup \big\{ (j,k,\ell) \in (\Z^2_+,\Z^2_-,\Z^2_+)\cup(\Z^2_-,\Z^2_-,\Z^2_-)\cup (\Z^2_+,\Z^2_+,\Z^2_-)  \;\big|\; \ell=j-k\big\}\\ 
  &\cup\big\{ (j,k,\ell) \in (\Z^2_-,\Z^2_+,\Z^2_+)\cup (\Z^2_-,\Z^2_-,\Z^2_-) \cup(\Z^2_+,\Z^2_+,\Z^2_-)\;|\; \ell=k-j\big\}  \\
    \Idx_-=&\big\{  (j,k,\ell) \in (\Z^2_-,\Z^2_+,\Z^2_+)\cup(\Z^2_+,\Z^2_-,\Z^2_+)\cup (\Z^2_-,\Z^2_-,\Z^2_-)  \;\big|\; \ell=j+k\big\}
\end{align*}
 Setting
\begin{align*}
  U^{t,\phi}(s,x)=\sum_{k\in \Z^2_0} \Uf_k(s) \base_k(x)\;, \mbox{ and }
  \phi(x)=\sum_{k\in \Z^2_0} \phi_k \base_k(x)\;, 
\end{align*}
we also have the backward equations
\begin{align}
  \label{eq:Ua}
  \begin{split}
    \frac{d}{ds}\Uf_\ell(s)=&\nu\abs{\ell}^2\Uf_\ell(s) +
    \sum_{(j,k,\ell)\in \Idx_+^*} c(j,\ell)\solf_j(t)
    \Uf_k(t)\\&\qquad\qquad\qquad- \sum_{(j,k,\ell)\in \Idx_-^*}
    c(j,\ell)\solf_j(t) \Uf_k(t) \quad s< t\;,
  \end{split}\\
  \Uf_\ell(t)=&\phi_\ell \ .\notag
  \end{align}
  where
\begin{align*}
  \Idx_+^*=& \big\{ (j,k,\ell) \in (\Z^2_+,\Z^2_+,\Z^2_-)\cup
  (\Z^2_-,\Z^2_+,\Z^2_+)\cup (\Z^2_+,\Z^2_-,\Z^2_+) \;|\; j+k+\ell=0 \big\}
  \\&
  \cup  \big\{ (j,k,\ell) \in
  (\Z^2_+,\Z^2_+,\Z^2_-)\cup(\Z^2_-,\Z^2_-,\Z^2_-)\cup (\Z^2_+,\Z^2_-,\Z^2_+)   \;|\; \ell=j-k\big\}\\
  &\cup \big\{ (j,k,\ell) \in (\Z^2_-,\Z^2_+,\Z^2_+)\cup (\Z^2_-,\Z^2_-,\Z^2_-)\cup (\Z^2_+,\Z^2_-,\Z^2_+)\;|\; \ell=j+k \big\}\\
  \Idx_-^*=& \big\{ (j,k,\ell) \in
  (\Z^2_-,\Z^2_+,\Z^2_+)\cup(\Z^2_+,\Z^2_+,\Z^2_-)\cup
  (\Z^2_-,\Z^2_-,\Z^2_-) \;|\; \ell=k-j\big\} \; .
\end{align*}  

  We now continue the proof of Proposition \ref{lem3.1}. Notice that
  the $\Uf_\ell$ are continuous in time for every $\ell \in
  \mathcal{Z}_0$, every $\phi \in \L^2$ and every realization of the
  stochastic forcing. Hence if $\Uf_\ell \equiv 0$ for some
  realization of noise, then $\phi_\ell=0$. (The notation $x\equiv 0$
  means $x(s)=0$, $s\in[0,t)$.)
  Thus to prove the lemma it would be sufficient to show that there
  existed a fixed set $\Omega_1$ with positive probability so that for
  any $\phi \in S_\infty$, if $\ip{\Mal(t) \phi}{\phi} = 0$ on
  $\Omega_1$ then $\Uf_\ell\equiv 0$ for all $\ell \in \Zf \cup
  \mathcal{Z}_\infty$. This will be proven inductively.
  
  The base case of the induction is given by Corollary
  \ref{l:baseCase}. In the present notation, it simply says that for
  any $\omega \in \Omega$ if $\ip{\Mal(t)\phi}{\phi} = 0$ for some
  $\phi \in \L^2$ then $\Uf_\ell \equiv 0$ for all $\ell \in \Zf$.
  In particular, for any $\omega \in \Omega$ if
  $\ip{\Mal(t)\phi}{\phi} = 0$ then $\Uf_\ell \equiv \Uf_{-\ell}
  \equiv 0$ for all $\ell \in \mathcal{Z}_0$.  The proof of Theorem
  \ref{thm:main} would then be complete if we show that there exists a
  single subset $\Omega_1 \subset \Omega$ of full measure so that if
  $\Uf_\ell \equiv \Uf_{-\ell} \equiv 0$ on $\Omega_1$ for all $\ell
  \in \mathcal{Z}_n$ then $\Uf_\ell \equiv \Uf_{-\ell} \equiv 0$ on
  $\Omega_1$ for all $\ell \in \mathcal{Z}_{n+1}$. This inductive step
  is given by the next lemma, which once proved completes the proof of
  Theorem \ref{thm:main}.
  
\begin{lemma}\label{l:recursion} There exists a fixed subset
  $\Omega_1$ of full measure so that for any $\phi \in \L^2$ and
  $\ell \in \Z_0^2$ if $\Uf_\ell\equiv 0$ and $\Uf_{-\ell} \equiv 0$
  on $\Omega_1$ then  for all $j \in \mathcal{Z}_0$ such that $j^\perp \cdot
  \ell\not=0$ and $\abs{j}\not=\abs{\ell}$
  \begin{equation*}
      \Uf_{\ell + j}\equiv\Uf_{-(\ell + j)}\equiv \Uf_{\ell-j}\equiv
     \Uf_{j-\ell}\equiv 0
  \end{equation*}
  on $\Omega_1$.
\end{lemma}
\bpf 

We begin with some simple observations which will be critical shortly.
Notice that from \eqref{eq:fourierAlpha} -- \eqref{eq:Ua} one sees that for
$\ell \in \Zf$, $\solf_\ell(s)$ has the form
\begin{align}
  \solf_\ell(s)&=\solf_\ell(0) + \int_0^s \gamma_\ell(r) dr + W_\ell(s) 
\end{align}
where the $\gamma_\ell$  are some stochastic processes
depending on the initial conditions and noise realizations.  Hence
these coordinates are the sum of a Brownian Motion and a part which
has finite first variation and is continuous in time for all $\omega
\in \Omega$.

Similarly, for $\ell \not\in \Zf$,
\begin{equation}
  \label{eq:XY}
  \begin{aligned}
    \solf_\ell(s)&=\solf_\ell(0) + \int_0^s \gamma_\ell(r) dr
  \end{aligned}
\end{equation}
and hence these coordinates are continuous and have finite first
variation in time for every $\omega \in \Omega$ as they are not directly forced. Similarly notice
that $\Uf_\ell$ is continuous and of finite first variation. In
particular, we emphasis that these properties of $\solf_\ell$ and
$\Uf_\ell$ hold on all of $\Omega$ for all $\ell \in \Z_0^2$ and $\phi
\in \L^2$.

Now, if $\Uf_\ell\equiv 0$ or $\Uf_{-\ell} \equiv 0$ then $
\frac{d}{ds}\Uf_\ell(s)\equiv 0$ or $\frac{d}{ds}\Uf_{-\ell}(s)\equiv 0$ 
respectively as
the coordinates are constant. Notice that from \eqref{eq:Ua} --
\eqref{eq:XY} these derivatives have the form
\begin{align*}
  X(s) + \sum_{k\in \Zf} Y_k(s) W_k(s)
\end{align*}
where the $X$ and $Y_k$ are continuous and bounded variation
processes. Also notice that they are \textbf{not adapted} to the past
of the $W_k$'s ! Nonetheless, it follows from Lemma \ref{lem3.2} in
the next section that if $\{X(\cdot), Y_k(\cdot): k \in \Zf\}$
are continuous and of bounded variation, then
\begin{equation}
\label{eq:XYWAssumption}
  X(s) +\sum_{k\in \Zf}Y_k(s) W_k(s)=0,\:\:0 \leq s \leq t,
\end{equation}
implies that
\begin{equation}\label{eq:XYWConclusion}
  Y_k(s)=0, \:\: \mbox{ for all $0 \leq s \leq t$ and  $k \in\Zf$}
\end{equation}
on a set $\Omega_1 \subset \Omega$, of full measure, which does not
depend on $\ell$, $k$ or $\phi$, and hence we can use a single
exceptional set for all of the steps in the induction. To summarize,
we have shown that there is a single fixed set $\Omega_1 \subset
\Omega$, of full measure, so that for any $\phi \in \L^2$ if
$\Uf_\ell\equiv \Uf_{-\ell}\equiv0$ on $\Omega_1$ then $Y_k\equiv 0$
for all $k \in\Zf$. We now identify the $Y_k$ to discover what
\eqref{eq:XYWConclusion} implies.

Define $|\ell|= \pm \ell$ depending on whether $\ell \in
\Z^2_\pm$ and $\sgn(\ell)=\pm 1$ depending on whether $\ell \in
\Z^2_\pm$. (Care should be taken not to confuse $\abs{\ell}$ which is
in $\R_+$ with $|\ell|$ which is in $\Z_+^2$.) 
Then from \eqref{eq:fourierAlpha} -- \eqref{eq:Ua}, we see that
for each $\ell \in \Z_0^2$
\begin{align*}
  \frac{d}{ds}\Uf_\ell(s) = X_\ell(s) & + \sum_{\substack{j \in \Zf \\
      (j,\ell) \in\Id_S}}
  c(j,\ell)\Big[\Uf_{-|\ell-j|}(s) +\sgn(\ell)
  \Uf_{-|\ell+j|}(s)\big]W_j(s)\\&+\sum_{\substack{j \in \Zf \\ 
      (j,\ell) \in \Id_A}} c(j,\ell) \Big[\Uf_{|\ell-j|}(s)
  -\sgn(\ell+j) \Uf_{|\ell+j|}(s) \Big]W_{j}(s),
\end{align*}
where $\Id_A= (\Z^2_+,\Z^2_-)\cup (\Z^2_-,\Z^2_+)$, $\Id_S=
(\Z^2_+,\Z^2_+)\cup (\Z^2_-,\Z^2_-)$, and $X_\ell(s)$ is a continuous
stochastic process with bounded variation. Hence by Lemma \ref{lem3.2}
we obtain that terms in brackets in the above equation are identically
zero.

Recall that by assumption $\frac{d}{ds}\Uf_\ell(s)=0$,
$\frac{d}{ds}\Uf_{-\ell}(s)=0$ and $\{j,-j\} \subset \Zf$.  Without
loss of generality, we assume that $\ell,j \in \Z_+$ since this can
always be achieved be renaming $\ell$ and $j$.  The preceding
reasoning using Lemma \ref{lem3.2} applied to $(j,\ell)$,
$(-j,-\ell)$, $(-j,\ell)$, and $(j,-\ell)$ implies respectively that
\begin{align*}
 c(j,\ell)\big[\Uf_{-|\ell-j|}(s) +
  \Uf_{-|\ell+j|}(s)\big]&=0\\
 c(j,\ell)\big[\Uf_{-|\ell-j|}(s) -
  \Uf_{-|\ell+j|}(s)\big]&=0\\
  c(j,\ell)\big[ \sgn(\ell-j)\Uf_{|\ell-j|}(s)-\Uf_{|\ell+j|}(s)\big]&=0\\
c(j,\ell)\big[ \sgn(\ell-j)\Uf_{|\ell-j|}(s)+\Uf_{|\ell+j|}(s)\big]&=0
\end{align*}
for all $s < t$ on a subset of $\Omega_1$ of full measure. Provided
that $j^\perp \cdot \ell \not=0$ and $\abs{j} \not=\abs{\ell}$, one has
that $c(\ell,j)\neq0$. Hence the left-hand-sides are linearly
independent and one concludes that $\Uf_{\ell -j}\equiv
\Uf_{\ell+j}\equiv \Uf_{j-\ell}\equiv \Uf_{-(\ell+j)}\equiv 0$ on a
subset of $\Omega_1$ of full measure. \epf

We now collect some of the information from the preceding proof for
later use.
\begin{proposition} \label{prop:formXYW} Let $U^{\phi,t}$ be the solution of
  \eqref{2.2} for any choice of terminal condition $\phi$ and terminal
  time $t$. Recall the definition of $S_n$ from \eqref{eq:S}.
  Let $\Pi_0$ be the projection onto $S_0$ and $\Pi^\perp_0$ its
  orthogonal complement.  Then for $s < t$
  \begin{align*}
    \frac{\partial\ }{\partial s}U^{\phi,t}(s)=X^\phi(s)+\sum_{j \in \Zf}
    Y_j^{\phi}(s) W_j(s)
  \end{align*}
  where
  \begin{align*}
    X^{\phi}(s)&=-\nu \Delta U^{\phi,t}(s) - B(\Pi^\perp_0
    \sol(s),U^{\phi,t}(s)) + C(U^{\phi,t}(s),\Pi^\perp_0 \sol(s))  \\&\qquad-
    B(R(s),U^{\phi,t}(s)) + C(U^{\phi,t}(s),R(s)),\\
    R(s)&=\Pi_0 \sol(0) + \int_0^s \nu\Delta \Pi_0 \sol(r) + \Pi_0
    B(\sol(r),\sol(r)) dr\\
    Y_j^{\phi}(s)&=-B(\base_j,U^{\phi,t}(s))+ C(U^{\phi,t}(s),\base_j)
  \end{align*}
  For all $\ell,j \in \Z^2_+$, we have 
  \begin{align*}
    \ip{ Y_j^{\phi}(s)}{\base_\ell}
    =&\pi^2c(j,\ell)\big[\Uf_{-|\ell-j|}(s) + \Uf_{-(\ell+j)}(s)\big]
    \\ 
    \ip{Y_{-j}^{\phi}(s)}{\base_{-\ell}} =&
   \pi^2c(j,\ell)\big[\Uf_{-|\ell-j|}(s) -\Uf_{-(\ell+j)}(s)\big]\\
    \ip{ Y_{-j}^{\phi}(s)}{\base_\ell} =&\pi^2c(j,\ell)\big[ \sgn(\ell-j)\Uf_{|\ell-j|}(s)-\Uf_{\ell+j}(s)\big]\\
    \ip{ Y_j^{\phi}(s)}{\base_{-\ell}}
    =&-\pi^2c(j,\ell)\big[ \sgn(\ell-j)\Uf_{|\ell-j|}(s)+\Uf_{\ell+j}(s)\big]\;.
  \end{align*}
  \end{proposition}
\section{A Quadratic Variation Lemma}
\label{sec:nonAdapted}
The following lemma is the main technical result used to prove the
existence of a density.  
\begin{lemma}\label{lem3.2}Let $\mathcal{A}$ be a collection of real
  valued stochastic processes such that there exists a fixed subset
  $\Omega_\mathcal{A} \subset \Omega$ of full measure such that on
  $\Omega_\mathcal{A}$ any element of $\mathcal{A}$ is continuous and
  has finite first variation.
  
  Fix a finite collection $\{W_1,\ldots,W_N\}$ of independent Wiener
  processes and a sequence of partitions $\{s_j^n\}_{j=0}^{m(n)}$, with
  $s_{j+1}^n -s_j^n \to 0$ as $n\to\infty$ and 
  \begin{equation*}
    0=s_1^n \leq  \cdots \leq s_{m(n)}^n=t \ .
  \end{equation*}
Then there exists a fix subset
  $\Omega' \subset \Omega_\mathcal{A}$ of full measure and a fixed subsequence of
  partitions $\{t_j^n\}_{j=0}^{k(n)}$ of  $\{s_j^n\}_{j=0}^{m(n)}$ so that if
\[Z(s)=X(s)+ \sum_{i=1}^N Y_i(s) W_i(s),\]
with $X, Y_1,\ldots, Y_N \in \mathcal{A}$ then on the set $\Omega'$
\begin{equation*}
  \sum_{j=1}^{k(n)} |Z( t_j^n)-Z(t_{j-1}^n)|^2  \longrightarrow \sum_{i=1}^N
  \int_0^t Y_i^2(s) ds  \text{  as $n \rightarrow \infty$. }
\end{equation*}
\end{lemma}
To prove this lemma we will invoke the following auxiliary results
whose proofs will be given after the proof of Lemma \ref{lem3.2}.
\begin{lemma}\label{lem3.3}
Let $\{W(s), 0 \leq s \leq t\}$ be a standard Brownian motion and
$\{s_j^n\}_{j=0}^{m(n)}$ a sequence of partitions as in Lemma \ref{lem3.2}. The sequence of measures
\begin{equation*}
  \left\{ \sum_{j=1}^{m(n)} (W(s_j^n)-W(s_{j-1}^n))^2 \delta_{s_{j-1}^n},\:\:
  n=1,2 \ldots\right\} 
\end{equation*}
converges weakly as $n \rightarrow \infty$ to the Lebesgue measure on
$[0,t]$, in probability.
\end{lemma}\begin{lemma}\label{lem3.4}
  Let $\{W(s),\:\:0 \leq s \leq t\}$ and $\{\bar{W}(s), \:\:0 \leq s
  \leq t \}$ be two mutually independent Brownian motions and
$\{s_j^n\}_{j=0}^{m(n)}$ a sequence of partitions as in Lemma \ref{lem3.2}. The
  sequence of signed measures
\begin{equation*}
\left\{ \sum_{j=1}^{m(n)}
  (W(s_j^n)-W(s_{j-1}^n))(\bar{W}(s_j^n)-\bar{W}(s_{j-1}^n))
  \delta_{s_{j-1}^n},\:\:n=1,2, \ldots  \right\} 
\end{equation*}
converges weakly to zero in probability as $n \rightarrow \infty$.
\end{lemma}

\bpf[Proof of Lemma \ref{lem3.2}.] In light of Lemma \ref{lem3.3} and
\ref{lem3.4}, since the collection of Brownian Motions is finite we
can select a single set of full measure $\Omega' \subset
\Omega_\mathcal{A}$ and a single subsequence of partitions
$\{t_j^n\}_{j=0}^{k(n)}$ such that the weak convergences given in
Lemma \ref{lem3.3} and \ref{lem3.4} hold on $\Omega'$. Furthermore we
can assume that on $\Omega'$ the $W_i(s)$ are continuous and have
quadratic variation $s$. For notational
brevity, we will write $t_j$ instead of $t_j^n$.  

Consider the quantity
$\sum_j|Z(t_j)-Z(t_{j-1})|^2$. If we express it in terms of the $X, Y$
and $W$'s we note that it contains four types of terms
\begin{equation*}
  \begin{split}
    &\sum_j|X(t_j)-X(t_{j-1})|^2\\
    &\sum_j|Y(t_j)W(t_j)-Y(t_{j-1})W(t_{j-1})|^2\\
    &\sum_j\left(X(t_j)-X(t_{j-1})\right)\left(Y(t_j)W(t_j)-Y(t_{j-1})W(t_{j-1})\right)\\
    &\sum_j\left(Y(t_j)W(t_j)-Y(t_{j-1})W(t_{j-1})\right)\left(\bar{Y}(t_j)\bar{W}(t_j)
      -\bar{Y}(t_{j-1})\bar{W}(t_{j-1})\right)\\
  \end{split}
\end{equation*}
where $W$ and $\bar{W}$ are mutually independent scalar Brownian
motions. The first and third terms are easily shown to tend to zero on
$\Omega'$ as
$n \rightarrow \infty$, since $X$ is of bounded variation, and $X$, $Y$
and $W$ are continuous.

Consider the second term :
\begin{equation*}
\begin{split}
&\sum_j|Y(t_j)W(t_j)-Y(t_{j-1})W(t_{j-1})|^2=\sum_j Y(t_{j-1})^2\left(W(t_j)-W(t_{j-1})\right)^2\\
&+\sum_jW^2(t_j)\left(Y(t_j)-Y(t_{j-1})\right)^2\\
&+2 \sum_j  W(t_j)Y(t_{j-1})\left(W(t_j)-W(t_{j-1})\right)\times \left(Y(t_j)-Y(t_{j-1})\right)
\end{split}
\end{equation*}
Again on $\Omega'$, the second and last terms above  tend to zero, and
\[\sum_j Y(t_{j-1})^2(W(t_j)-W(t_{j-1}))^2 \rightarrow \int_0^t Y(s)^2 ds\]
on $\Omega'$
 by the convergence given in Lemma \ref{lem3.3}.
Finally
\begin{equation*}
  \begin{split}
    &\sum_j(Y(t_j)W(t_j)-Y(t_{j-1})W(t_{j-1}))(\bar{Y}(t_j)\bar{W}(t_j)-\bar{Y}(t_{j-1})\bar{W}(t_{j-1}))\\
    &=\sum_j \left[Y(t_{j-1})(W(t_j)-W(t_{j-1}))+(Y(t_j)-Y(t_{j-1})) W(t_j)\right]\\
    &\times\left[\bar{Y}(t_{j-1})(\bar{W}(t_j)-\bar{W}(t_{j-1}))+
      (\bar{Y}(t_j)-\bar{Y}(t_{j-1}))\bar{W}(t_j)\right]\\
    &=\sum_j Y(t_{j-1})\bar{Y}(t_{j-1})\left(W(t_j)-
      W(t_{j-1})\right)\left(\bar{W}(t_j)- \bar{W}(t_{j-1})\right)
    + \varepsilon_n,
  \end{split}
\end{equation*}
where $\varepsilon_n \rightarrow 0$ a.s., as $n \rightarrow \infty$.
Again by Lemma \ref{lem3.4}, the sum tends to zero on $\Omega'$.\epf

\bpf[Proof of Lemma \ref{lem3.3}.] For any measure, the fact that
$\mu_n \Rightarrow \mu$ follows from $\mu_n([0,s])\rightarrow
\mu([0,s])$ for all $s \in [0,t]$, $s$ rational. From any subsequence
of the given sequence, one can extract a further subsequence such that
$\mu_n([0,s]) \rightarrow \mu([0,s])$, for all $s$ rational, $0 \leq s
\leq t$, a.s. Hence along that subsequence $\mu_n \Rightarrow \mu$
a.s., hence the whole sequence converges weakly in probability.  \epf
\bpf[Proof of Lemma \ref{lem3.4}.] Write $\Delta_j W$ for $W(s_j^n)
-W(s_{j-1}^n)$ and $\Delta_j \bar{W}$ for
$\bar{W}(s_j^n)-\bar{W}(s_{j-1}^n)$. Note that for all $0 \leq r < s
\leq t$, as $n \rightarrow \infty$,
\[\sum_{r < s_{j}^n \leq s} \Delta_j  W \Delta_j \bar{W} \rightarrow 0
\mbox{  in probability, as } n\rightarrow\infty.\] 
Consequently if $f$ is a step function, 
\[ \sum_{j=1}^n \Delta_j W \Delta_j \bar{W} f(s_{j-1}^n) \rightarrow 0
\mbox{  in probability, as } n \rightarrow \infty.\] 
 Moreover for any two functions $f$ and $g$
\begin{equation*}
\begin{split}
  &\Big|\sum_{j=1}^n \Delta_jW \Delta_j
  \bar{W}\Bigl(f(s_{j-1}^n)-g(s_{j-1}^n)\Bigr)\Big|\\ 
  & \leq \sup_{0 \leq s \leq
    t}\big|f(s)-g(s)\big|\Bigl(\sum_{j=1}^n(\Delta_jW)^2\Bigr)^{\hf} \Bigl( 
    \sum_{j-=1}^n(\Delta_j \bar{W})^2 \Bigr)^{\hf},
\end{split}
\end{equation*}
and the right hand side tends to
\begin{equation*}
t \times \sup_{0 \leq s \leq t} |f(s) -g(s)|
\end{equation*}
in probability, as $n \rightarrow \infty$.

Let now  $f$ be a continuous function, and $g$ be a step function.
Choose
\[\displaystyle\delta=2 t \sup_{0 \leq s \leq t}|f(s) -g(s)|\]
We have
\begin{equation*}
\begin{split}
  {\P}\Bigl\{\Big|\sum f(s_{j-1}^n) \Delta_j W \Delta_j\bar{W}\Big| >
    \delta\Bigr\}
  &\leq \P\Bigl\{\Big|\sum_{j=1}^n g(s_{j-1}^n) \Delta_j W \Delta_j
    \bar{W}\Big| > \delta/3\Bigr\}\\ 
  +{\P}&\Bigl\{\Big|\sum_{j=1}^n\left(f(s_{j-1}^n)-g(s_{j-1}^n)\right)\Delta_j
    W \Delta_j \bar{W}\Big| > 2 \delta/3 \Bigr\},
\end{split}
\end{equation*}
and it follows from the above arguments that the latter tends to zero
as $n \rightarrow \infty$. Since $\delta$ can be made arbitrarily
small by an appropriate choice of the step function $g$, the lemma is
proved. \epf

\section{Relation to Brackets of vector fields}
\label{sec:horomander}
We now sketch another possible proof of our  Theorem
\ref{thm:main}, which brings in explicitly the brackets of certain
vector fields.  A vector field over the space $\L^2$ is a mapping from
a dense subset of $\L^2$ into itself.  We begin by rewriting
\eqref{eq:vorticity} as
\begin{equation*}
  \frac{\partial\sol}{\partial t}(t)=F_0(\sol(t)) + \sum_{i=1}^N F_i
  \frac{\partial W_i}{\partial t}(t) .
\end{equation*}
The diffusion  vector fields in our case are constant vector fields
defined by
\begin{equation*}
  F_i=\base_{k_i},\quad 1\le i\le N,  
\end{equation*}
where  $N$ is
the cardinality of $\Zf$ and $\{k_1,\cdots,k_N\}$ is any ordering of the set $\Zf$.  Similarly the drift vector field is denoted
by $F_0(\sol)=\nu \Delta \sol - B(\sol,\sol)$. In this notation,
\eqref{prop2.2}, becomes
\begin{equation*}
  \left\{
    \begin{aligned}
     &\frac{\partial U^{t,\phi}}{\partial s}(s) + (\nabla_\sol F_0)^*(\sol(s))U^{t,\phi}(s) =0, \qquad
      \hfill  0 \leq s \leq t;\\
      &U^{t,\phi}(t)=\phi,
    \end{aligned}
  \right.
\end{equation*}
where $\nabla_\sol F_0$ is the Fr\'echet derivative of $F_0$ in the $\L^2$
topology and $ (\nabla_x F_0)^*$ is its $\L^2$-adjoint. If it is well
defined, we define the bracket $[F,G]$ between two $\L^2$ vector fields
$F$ and $G$ as $[F,G]=(\nabla_\sol F)G- (\nabla_\sol G)F$. (Part of
being well defined is that the range of $G$ and $F$ are contained in
the domain of $\nabla_\sol F$ and $\nabla_\sol G$ respectively.)

The argument in this alternate proof is based on the two next results.
\begin{lemma}
  Let $G$ be a vector field on $\L^2$ which is twice Fr\'echet
  differentiable in the $\L^2$ topology and such that $[F_i,G],
  i=0\dots N$, are vector fields on $\L^2$. Then we have
\begin{align*}
  &\ip{U^{t,\phi}(s)}{G(\sol(s))}-\ip{U^{t,\phi}(0)}{G(\sol(0))}\\
  &= \int_0^s\ip{U^{t,\phi}(r)}{[F_0,G](\sol(r))}dr
  +\sum_{i=1}^N\int_0^s\ip{U^{t,\phi}(r)}{[F_i,G](\sol(r))}\circ dW^i_r\\
  \intertext{where $\circ$ means that it is a Stratonovich (anticipating) integral;
    in It\^o--Skorohod language, it takes the form}
  &=\int_0^s\ip{U^{t,\phi}(r)}{[F_0,G](\sol(r))}dr
  +\sum_{i=1}^N\int_0^s\ip{U^{t,\phi}(r)}{[F_i,G](\sol(r))} dW^i_r\\
  &+\sum_{i=1}^N\int_0^s\left[\frac{1}{2}\ip{ \nabla_\sol^2
      G(\sol(r))(F_i,F_i)}{U^{t,\phi}(r)}
    +\ip{[F_i,G](\sol(r))}{D^i_rU^{t,\phi}(r)}\right]dr
\end{align*}
\end{lemma} 

\bpf The formula in Skorohod language follows from Theorem 6.1 in 
\cite{b:NualartPardoux88}, via an easy finite dimensional 
approximation. 
Its translation in the Stratonovich form follows 
from
Theorem 7.3 in the same paper (see also Theorem 3.1.1 in \cite{b:Nualart95}). 
\epf

We can now prove the following:
\begin{proposition}
Let $\Omega_0\subset\Omega$. 
Under the same assumptions on $G$ as in the above Lemma,
\begin{align*}
  \ip{U^{t,\phi}(s)}{G(\sol(s))}&\equiv 0 \text{ on the set } \Omega_0\\
  \intertext{implies that} \ip{U^{t,\phi}(s)}{[F_i,G](\sol(s))}&\equiv
  0 \text{ a.s. on the set } \Omega_0,\, i=0,1,\ldots,N.
\end{align*}
\end{proposition}
\bpf 
The assumption implies that the quadratic variation on $[0,t]$ of the
process \\
$\{\ip{U^{t,\phi}(s)}{G(\sol(s))}\}$ vanishes almost surely on $\Omega_0$. Then from \cite{b:Nualart95},
Theorem 3.2.1, for $i=1,\ldots,N$,
\begin{equation*}
  \int_0^t \ip{U^{t,\phi}(s)}{[F_i,G](\sol(s))}^2 ds = 0  
\end{equation*}
a. s. on the set $\Omega_0$, i.e.
\begin{equation*}
  \ip{U^{t,\phi}(s)}{[F_i,G](\sol(s))}\equiv 0 \text{ a. s. on the set
  } \Omega_0,\,  i=1,\ldots,N.  \
\end{equation*}
This implies that for $1\le i\le N$,
\begin{equation*}
  \int_0^s\ip{U^{t,\phi}(r)}{[F_i,G](\sol(r))}\circ dW^i_r\equiv 0\text{
    a. s. on the set } \Omega_0    
\end{equation*}
(see Definition 3.1.1 in \cite{b:Nualart95}), from which it follows
(see the previous Lemma) that
\begin{equation*}
  \ip{U^{t,\phi}(s)}{[F_0,G](\sol(s))}\equiv 0 \text{ a. s. on the set } \Omega_0.   
\end{equation*}
\epf

Now call $\mathcal{L}$ all well defined $\L^2$ vector fields in the
ideal generated by the vector fields $F_1,\ldots,F_N$ in the Lie
algebra generated by $F_0,F_1,\ldots,F_N$.  In other words, at each
$u\in\L^2$, $\mathcal{L}(u)$ consists of $F_1,\ldots,F_N$, and all
brackets
\begin{equation}
  \label{eq:brackets}
  [F_{i_n},[F_{i_{n-1}},\ldots[F_{i_2},F_{i_1}]\ldots](u),  
\end{equation}
which are well defined vector fields on $\L^2$ where $1\le i_1\le N$,
and for $j>1$, $0\le i_j\le N$.  Iterating the argument in the
Proposition, we deduce the following result.
\begin{corollary}
Given $\phi\in \L^2$, let
\begin{equation*}
  \Omega_0:=\{\ip{\Mal(t)\phi}{\phi}=0\}.  
\end{equation*}
If $\P(\Omega_0)>0$, then 
\begin{equation*}
  \ip{\phi}{G(\sol(s))}\equiv 0,\quad\forall G\in\mathcal{L},\, \text{
    a. s. on }\Omega_0. 
\end{equation*}
In particular, $\phi$ is orthogonal to all constant vector fields in $\mathcal{L}$.
\end{corollary}
In the case of the stochastic Navier Stokes equations given in
\eqref{eq:vorticity}, all of the brackets given in \eqref{eq:brackets}
are well defined if the $e_k$ used to define the forcing have
exponentially decaying Fourier components. This follows from the fact
that all of the bracket of the form \eqref{eq:brackets} contain
differential operators with polynomial symbols and the fact that, with
this type of forcing, on any finite time interval $[0,T]$ there exists
a positive random variable $\gamma$ so that $\sup_{[0,T]}
\|e^{\gamma|\nabla|}\sol(s)\|^2 < \infty$ almost surely. Here
$\|e^{\gamma|\nabla|}\sol\|^2= \sum_k e^{2\gamma|k|}|\sol_k|^2$
where $\sol(t,x)=\sum_k \sol_k(t) e^{i k \cdot x}$. See for example
\cite{b:Mattingly02a} for a stronger version of this result or
\cite{b:MattinglySinai98,b:Mattingly98b} for simpler versions.

Furthermore in \cite{b:EMattingly00} it was implicitly shown by the
construction used that the span of the constant vector fields contains
$S_\infty$. Thus, under the same conditions as before we see that the
law of arbitrary finite dimensional projections of $\sol(t)$ have a
density with respect to Lebesgue measure.
\section{Smoothness}\label{sec:smooth}
In the preceding sections, we proved the existence of a density. We now
address the smoothness of the density. While the former simply
required that the projected Malliavin matrix be invertible, the proof of
smoothness requires control on the norm of the inverse of the
projected Malliavin matrix  together with ``smoothness in the
Malliavin sense.'' The following is the main result of
this section; however, it rests heavily on the general results proven
in Section \ref{sec:noSmallEVs}, as well as some technical results
from the appendices.

\begin{theorem}\label{thm:mainSmooth}
 Let $S$ be any finite dimensional subspace of $S_\infty$ and $\Pi$
  the orthogonal projection in $\L^2$ onto $S$. For any $t>0$, the law
  of $\Pi w_t$ has a $C^\infty$ density with respect to the Lebesgue
  measure on  $S$.
\end{theorem}
\bpf We use Corollary 2.1.2 of \cite{b:Nualart95}. Lemma
\ref{l:Dinfty} from the appendix  establishes condition \textit{(i)}
from that corollary while condition \textit{(ii)} of the same
corollary follows from the
next theorem.\epf

The following is a quantitative version of Proposition \ref{lem3.1}. It
gives a quantitative control of the smallest eigenvalue of a finite
dimensional projection of the Malliavin matrix.
\begin{theorem} \label{l:onePhi} Let $\Pi$ be the orthogonal projection of
  $\L^2$ onto a
  finite dimensional subspace of $S_\infty$. For any $T>0$, $\eta >0$,
  $p \geq 1$, and $K > 0$ there exists a constant
  $c=c(\nu,\eta,p,|\Zf|,T,K,\Pi)$ and 
  $\epsilon_0=\epsilon_0(\nu,K,|\Zf|,T,\Pi)$ so that for all $\epsilon \in
  (0,\epsilon_0]$,
  \begin{align*}
    \P\Bigl( \inf_{\phi \in S(K,\Pi)}
    \ip{\Mal(T)\phi}{\phi} < \epsilon \Bigr) \leq c\exp(\eta \|\sol(0)\|^2) \epsilon^p
  \end{align*}
where $S(K,\Pi)=\{ \phi \in S_\infty: \| \phi \|_1\leq 1, \| \Pi
\phi\| \geq K \}$.
\end{theorem}
\begin{remark}
  Notice that this lemma implies that the 
  eigenvectors with ``small'' eigenvalues have small projections in the
  ``lower'' modes. The definition of ``lower'' modes depends on
  the definition of  ``small'' eigenvalues. This separation between
  the eigenvectors with small eigenvalues and the low modes is one of
  the keys to the ergodic results proved in \cite{b:HairerMattingly04}.
\end{remark}
\begin{remark}
  Also notice that there is a mismatch in the topology in Theorem
  \ref{l:onePhi} in that the test functions are bounded in $\H^1$ but
  the innerproduct is in $\L^2$. This can likely be rectified since
  the backward adjoint linearized flow $\bar J_{s,T}^*$ maps $\L^2$
  into $\H^1$ for any $s \leq T$, it is possible to obtain estimates
  on $\P(\langle \phi, \Mal(T) \phi\rangle < \epsilon )$ for $\phi \in
  \L^2$. Just as in the proof of Theorem \ref{l:onePhi}, where we 
  exclude a small neighborhood of time zero to allow $w_t$ to
  regularize, we could exclude the $s$ in a small neighborhood of the
  terminal $T$ to allow $U^{T,\phi}(s)=J_{s,T}^*\phi$ to regularize.
\end{remark}
\bpf 
Recall the definition of $S_n$ and $\mathcal{Z}_n$ from
\eqref{eq:S} and let $\Pi_n$ be the orthonormal projection onto $S_n$.
Since $S(K,\Pi) \subset S_\infty$ and $\Pi$ projects onto a finite
dimensional subspace of $S_\infty$, for $n$ sufficiently large $\|\Pi_n \phi\| >
\frac12 K$ for all $\phi \in S(K,\Pi)$. Fix such an $n$.

We now construct a basis of $S_n$ compatible with the structure of
$\mathcal{Z}_k$, $k\leq n$.  Fixing any ordering of $\Zf$, set $\{f_i:
i=1,\dots,N=|\Zf|\}=\Zf$. Clearly $\{f_i\}_{i=1}^N$ is a basis for
$S_0$. Set $J_0=N$. By the construction of $S_n$ it is clear that
$S_n\setminus S_{n-1}$ is equal to the $\mbox{span}(\base_k: k\in
\mathcal{Z}_n')$ where $\mathcal{Z}_n'\eqdef\mathcal{Z}_n
\cap_{j=0}^{n-1} \mathcal{Z}_{j}^c\cap\mathcal{Z}_*^c$.  For $n\geq
1$, set $J_n=J_{n-1}+|\mathcal{Z}_n'|$ and $\{f_i :
i=J_{n-1}+1,\dots,J_n\}=\{ \base_k : k\in\mathcal{Z}_n'\}$, again fixing
an arbitrary ordering of the righthand side. Clearly $\{
f_i\}_{i=1}^{J_n}$ chosen in this way is an orthogonal basis for
$S_n$.

Fix some $t_0 \in (0,T)$. Recall that by Corollary \ref{cor2.4}
\begin{equation*}
     \P\Bigl( \inf_{\phi \in S(K,\Pi)}
    \ip{\Mal(T)\phi}{\phi} < \epsilon \Bigr)=\P\Bigl(\inf_{\phi \in
      S(K,\Pi)}\sum_{k \in \Zf} \int_0^T \ip{U^{T,\phi}(s)}{\base_k}^2
    ds < \epsilon \Bigr)\;.
\end{equation*}
Let $X^\phi$ and $Y_j^{\phi}$ be as in
proposition \ref{prop:formXYW}. For $j=1,\cdots,J_n$, define
\begin{equation*}
\chi_j^\phi(t)=-\ip{X^\phi(T-t)+\sum_{k \in \Zf} Y_k^{\phi}(T-t)
  W_k(T)}{f_j} \;,   
\end{equation*}
$\Upsilon_{j,k}^{\phi}(t)=\ip{Y_k^{\phi}(T-t)}{f_j}$, 
$G_j^\phi(t)=\ip{U^{T,\phi}(T-t)}{f_j}$, and $\bar
W_k(t) = W_k(T)-W_k(T-t)$.
Notice that for $t \in [0,T]$ and  $j=1,\cdots,J_n$  
\begin{align*}
  G_j^\phi(t)=\ip{\phi}{f_j} + \int_0^t \Big[\chi_j^\phi(s) +\sum_{k
  \in \Zf}
  \Upsilon_{j,k}^{\phi}(s)   \bar W_k(s)  \Big] ds \ .
\end{align*}
Furthermore in light of the observations in the proof of Lemma
\ref{l:recursion}, we see that this sequence of equations satisfies
the assumptions of the next subsection
with $J_0$ as defined above and $J=J_n$. Next recall that
$U^{t,\phi}(s,x)=\sum_{k\in \Z^2_0} \Uf_k(s) \base_k(x)$. Combining this,
the last equalities in Proposition \ref{prop:formXYW}, and the
argument already used at the end of the proof of Lemma
\ref{l:recursion}, we see that each $G_j^\phi$ is a linear combination
of the $\{\Upsilon_{i,k}^{\phi} : i<j, k \in \Zf\}$ with coefficients
which are constant in time.

Below $\N{\cdot}_{1,[a,b]}$ and $\|\cdot\|_{\infty,[a,b]}$
respectively denote the Lipschitz and $L^\infty$ norm on $[a,b]$, see
the beginning of section \ref{sec:noSmallEVs} for the precise definitions.

Fix a $t\in (0,T)$ and set $T_1=T-t$. Bounds, uniform for $\phi\in
S(K,\Pi)$, on the $p$-th moments of the $L^\infty$--norm and Lipschitz
constants of $\chi_j^\phi$, and $\Upsilon_{j,k}^{\phi}$ over the
interval $[0,T_1]$ are given by Lemma \ref{l:XYMomentsForSNS}. (Recall
that these processes have been time reversed.) Hence given any $p\geq
1$, $q > 0$ and $\eta>0$, there exists a
$c=c(\eta,q,p,t,\nu,\mathcal{E}_1,T)> 0$ so that for any $\epsilon >0$
if one defines
\begin{align}\label{eq:Oflat}
  \Omega_\flat(\epsilon,q) =\bigcap_{\phi \in S(K,\Pi)} \left\{
    \sup_{k,j}( \N{\chi_j^\phi}_{1,[0,T_1]},
    \N{\Upsilon_{j,k}^{\phi}}_{1,[0,T_1]}) \leq
    \epsilon^{-q}\right\}
\end{align}
then the estimate $\P(\Omega_\flat(\epsilon,q)^c) \leq c \exp( \eta \|\sol(0)\|^2 )
\epsilon^p$ holds.

Next Corollary \ref{c:LpLadder} and Proposition \ref{prop:orderEpsp}
 state that there exist $q=q(|S|,N)$
and $\epsilon_0=\epsilon_0(T,\Zf|,|S|)$ so that for all $\epsilon \in
(0,\epsilon_0]$ there is a $\Ospp(\epsilon)$ so that for all $\phi \in
S(K,\Pi)$ one has
\begin{equation}\label{eq:OmegaSmall}
\left\{  \sup_{i=1,\cdots,N} \int_0^{T_1}
      |G_i^\phi(s)|^2 ds < \epsilon ; \sup_{i=1,\cdots,J_n}
      \sup_{s\in[0,T_1]}|G_i^\phi(s)| > \epsilon^q \right\} \cap
    \Omega_\flat(\epsilon,q) \subset \Ospp(\epsilon)
\end{equation}
and $\P(\Ospp(\epsilon)) \leq c\epsilon^p$ for all $p \geq 1$ and
$\eta>0$ with a $c=c(T,|\Zf|,|S|,p,\eta,\nu)$. Notice that because of the
uniformity in \eqref{eq:Oflat}, $\Ospp(\epsilon)$ does not depend on
the sequence of $G$'s. Since $
\sup_{k\in\Zf}\int_t^T\ip{U^{T,\phi}(s)}{e_k}^2ds = \sup_{i=1,\cdots,N}\int_0^{T_1} |G_i^\phi(s)|^2 ds$
and $ \sup_{s\in[0,T_1]}|G_i^\phi(s)| =
\|\ip{U^{T,\phi}}{f_i}\|_{\infty,[t,T]}$, the inclusion given in
\eqref{eq:OmegaSmall} becomes
\begin{align}\label{eq:supG1}
  \left\{ \sup_{k\in\Zf} \int_t^T \ip{U^{T,\phi}(s)}{\base_k}^2 ds <
    \epsilon ; \sup_{i=1,\cdots,J_n}
    \|\ip{U^{T,\phi}}{f_i}\|_{\infty,[t,T]} > \epsilon^q \right\}
  \subset \Ospp(\epsilon) \cup \Omega_\flat(\epsilon,q)^c
\end{align}
for all $\epsilon \in (0,\epsilon_0]$. Since
$\ip{U^{T,\phi}(T)}{f_i}=\ip{\phi}{f_i}$, by the choice of the
subspace $S_n$ one has $ \sup_{i} \|\ip{U^{T,\phi}}{f_i}\|_{\infty,[t,T]}
\geq \frac{K}{2\sqrt{J_n}}$. Thus for $\epsilon \in (0,\epsilon_0\wedge
(\frac{K}{2\sqrt{J_n}})^{\frac1q}]$ one has $\epsilon^q \leq \frac{K}{2\sqrt{J_n}}$
which transforms \eqref{eq:supG1} into
  \begin{align*}
    \Bigl\{ \sum_{k\in\Zf} \int_t^T \ip{U^{T,\phi}(s)}{\base_k}^2 ds <
      \epsilon\Bigr\}\subset   \Ospp(\epsilon) \cup
      \Omega_\flat(\epsilon,q)^c \ .
  \end{align*}
  As $\phi$ was an arbitrary direction in $S(K,\Pi)$ and
  $\Ospp(\epsilon)$ and $\Omega_\flat(\epsilon,q)$ are independent
  of $\phi$, we have
  \begin{multline*}
    \Bigl\{\inf_{\phi \in S(K,\Pi)} \sum_{k\in\Zf} \int_t^T
    \ip{U^{T,\phi}(s)}{\base_k}^2 ds < \epsilon\Bigr\}= \\\bigcup_{\phi
      \in S(K,\Pi)}\Bigl\{ \sum_{k\in\Zf} \int_t^T
    \ip{U^{T,\phi}(s)}{\base_k}^2 ds < \epsilon\Bigr\}\subset
    \Ospp(\epsilon) \cup \Omega_\flat(\epsilon,q)^c \ .
  \end{multline*}
In summary we have shown that for any $p\geq 1$ and $\eta >0$, there
exists $q >0$ so that the above inclusion holds and a $c > 0$ so that
 $\P(\Ospp(\epsilon)) + \P( \Omega_\flat(\epsilon,q)^c) \leq c
  \exp(\eta \|\sol(0)\|^2) \epsilon^p$.  \epf



\section{Controlling the Chance of Being Small}
\label{sec:noSmallEVs}
This section contains the main estimate used to control
the chance of certain processes being small when their quadratic
variation is large. The estimates of this section are simply
quantitative versions of the results of Section
\ref{sec:nonAdapted}. There are also the analog of results used in the
standard Malliavin calculus as applied to finite dimensional SDEs. There
the estimates were developed by Stroock \cite{b:Stroock83} and Norris
\cite{b:Norris86}. Here we do not have adapted processes. Instead, we
exploit the smoothness in time to obtain estimates.

For the entirety of this section, we fix a time $T$ and consider only
the interval of time $[0,T]$.
For any real-valued function of time $f$, define the $\alpha$-H\"older constant over the
time interval $[0,T]$ by 
\begin{align}\label{eq:DefLip}
  \hol_\alpha(f)= \sup_{\substack{s,r \in[0,T]\\0< |s-r|\leq 1}}
  \frac{|f(s)-f(r)|}{|s-r|^\alpha}
\end{align}
and the $L^\infty$ norm by 
\begin{align}\label{eq:defInf}
\|f\|_{\infty}=\sup_{s \in[0,T]}|f(s)| .  
\end{align}
We also define $\N{f}_{\alpha}=\max( \|f\|_{\infty}, \hol_\alpha(f))$.
At times we will also need versions of the above norms over shorter
intervals of time. For $[a,b] \subset [0,T]$ we will write
$\hol_{\alpha,[a,b]}(f), \|f\|_{\infty,[a,b]},$ and
$\N{f}_{\alpha,[a,b]}$ for the norms with the same definitions as
above except that the supremum over $[0,T]$ is replaced with a
supremum over $[a,b]$. We also extend the definitions of the Lipschitz
constant in time $\hol_\alpha(f)$, to functions of time talking values
in $\L^2$ by replacing the absolute value in the definitions given in
\eqref{eq:DefLip} and  \eqref{eq:defInf} by the norm on $\L^2$. Similarly we extend the
definition $\hol_{\alpha,[a,b]}(f)$, $\|f\|_\infty$,
$\|f\|_{\infty,[a,b]}$ to functions of time taking values in $\L^2$.

\subsection{A Ladder of Estimates}\label{sec:ladderDef}
For each $1 \leq j \leq J$, let $\{W^{(j)}_i : i=1\dots N\}$ be a collection of mutually
independent standard Wiener processes with $W_i(0)=0$ 
defined on a probability
space $(\Omega,\mathcal{F},\P)$.

We say that the collection of processes $\mathcal{G}=\{
G^{(j)}(t,\omega) : 1 \leq j \leq J\}$ forms a ladder of order $J$ with
base size $J_0$ if first $1 \leq J_0 <J$ and
\begin{align*}
  G^{(j)}(t,\omega)&=G_0^{(j)}+\int_0^t H^{(j)}(s,\omega) ds \\
  H^{(j)}(s,\omega) &= X^{(j)}(s,\omega) + \sum_{i=1}^N Y_i^{(j)}(s)
  W_i^{(j)}(s,\omega)
\end{align*}
where $j=1,\dots,J< \infty$ and $\omega \in \Omega$. 

Second, we require that for $j$ greater than $J_0$, the $G^{(j)}$ are
determined by the functions at the previous levels. More
precisely, for each $j$ with $j > J_0$ there exists an integer
$K=K(j)$, a collection $\{ g_k(t) : k=1,\dots,K\}$ of bounded,
deterministic functions of time, and a collection $\{ f_k(t,\omega) :
k=1,\dots,K\}$ of stochastic process with
 \begin{equation}\label{eq:earlyInLadder}
   f_k \in \{ Y_i^{(l)},   X^{(l)},  G^{(n)} : 1\leq i \leq N,1 \leq l
   \leq j-1 , 1\leq n \leq j-1 \}  
 \end{equation}
so that
\begin{align*}
  G^{(j)}(t,\omega) = \sum_{k=1}^K g_k(t)f_k(t,\omega)  \mbox{
  almost surely.}   
\end{align*}
This assumption can be restated by saying that, for $j >J_0$,
$G^{(j)}$  
must be at each moment of time in the span of the preceding $X$, $Y$,
and $G$. And furthermore, the coefficients in
the linear combination producing $G^{(j)}$ must be uniformly bounded
on $[0,T]$.

It is important to remark that we do \textrm{not} assume that the
$Y_i^{(j)}$ or $X^{(j)}$ are adapted to the Wiener processes. Typical
assumptions regarding adaptedness will be replaced with assumptions on
the regularity of the processes in time.

The goal of this section is to prove that under certain assumptions,
if the first $J_0$ of the  $G^{(j)}$ are small in some sense then all of the
$X$, $Y$, and remaining $G$ are also small with high probability. The
ladder structure connects the $j$-th level with the other levels.

Fix a time $T >0$. For any choice of the positive parameter $\Delta$
define $\delta=\Delta^\frac53$. For $k=0,1,\dots$, define $t_k=k\Delta
\wedge T$. For each fixed $k$, define $s_\ell(k)=(t_k + \ell \delta)
\wedge t_{k+1}$ for $\ell=0,1,\dots$. Set
$\delta_\ell^k=s_\ell(k)-s_{\ell-1}(k)$ and $\delta_\ell^k f =
f(s_\ell(k)) - f(s_{\ell-1}(k))$.  Lastly define $m=\inf\{ k : t_k=T\}$
and $M(k)=\inf\{ \ell : s_\ell(k)=t_{k+1} \}$.  Notice that
$\Delta^{-\frac23}=\frac{\Delta}{\delta}\leq M(k) \leq
\frac{\Delta}{\delta}+1= \Delta^{-\frac23}+1$ for all $k$ and $m \leq T
\Delta^{-1}+1$.

Define the following subsets of the
probability space $\Omega$:
\begin{align*}
  \Omega_a(\Delta)&=\Bigg\{ \inf_{0\leq k\leq m} \inf_{1\leq i\leq N}
  \frac1{M(k)} \sum_{\ell=1}^{M(k)} \frac{(\delta_\ell^k
    W_i)^2}{\delta_\ell^k}
  \leq \frac12\Bigg\}\\
  \Omega_b(\Delta)&=\Bigg\{ \sup_{0\leq k\leq m} \sup_{\substack{(i,j)
      \in \{1,\dots,N\}^2 \\ i\neq j}} \frac1{M(k)} \Big|
  \sum_{\ell=1}^{M(k)} \frac{(\delta_\ell^k W_i)(\delta_\ell^k
    W_j)}{\delta_\ell^k}\Big|
  \geq \frac{\Delta^{\frac{3}{14}}}{3N^2}\Bigg\}\\
  \Omega_c(\Delta)&=\Bigg\{ \sup_{1\leq i\leq N} \N{W_i}_\frac14 >
  \Delta^{-\frac1{28}} \Bigg\}\;,
\end{align*}
$\Os(\epsilon)= \Omega_a(\epsilon^\frac{14}{75}) \cup
\Omega_b(\epsilon^\frac{14}{75}) \cup
\Omega_c(\epsilon^\frac{14}{75})$ and finally
\begin{align*}
\Ospp(\epsilon)=\bigcup_{j=1}^J \Os(\epsilon^{(\frac{1}{152})^j})\;.  
\end{align*}
The following bound follows readily from Corollary \ref{c:ab} and
Lemma \ref{l:c} in section \ref{sec:SmallABC}.
\begin{proposition} \label{prop:orderEpsp}For any $p \geq 1$, there exists a constant
  $c=c(T,J,N,p)$, in particular independent of $\epsilon$, such that 
  \begin{equation*}
    \P\Bigl(\Ospp(\epsilon)\Bigr) \leq c \epsilon^p \;.
  \end{equation*}
\end{proposition}
The next key result in this section is the following proposition which
shows why the previous estimate is important.
\begin{proposition} \label{l:ladderEstimateLocalized}
  Fix a positive integer $J$ and for $q > 0$ define
 \begin{align*} 
   \Omega_{q*}(\mathcal{G},\epsilon)&=\Big\{
   \sup_{\substack{1\leq j\leq J\\
       1\leq i\leq N}}\big(\N{X^{(j)}}_{1},\N{Y_i^{(j)}}_{1}\big)\leq
   \epsilon^{-q}\Big\}\;,\\
   \Omega_q(\mathcal{G},\epsilon)&=\Biggl\{ \sup_{1\leq j\leq J_0}
   \|G^{(j)}\|_{\infty} < \epsilon \ ;\sup_{\substack{1\leq j\leq J\\
       1\leq i\leq N}}\big(\|X^{(j)}\|_{\infty},
   \|Y_i^{(j)}\|_{\infty}, \|G^{(j)}\|_\infty\big)>\epsilon^q \Biggr\}
   \; .
  \end{align*} 
  Then there exists positive constants
  $q_0=q_0(J)$  and $\epsilon_0=\epsilon_0(T,J,N)$ so that for
  any $\epsilon \in (0,\epsilon_0]$, $q \in (0,q_0]$ and 
  ladder $\mathcal{G}=\{ G^{(j)} : 1\leq j \leq J \}$ of order $J$
  with base size $J_0$ less than $J$, 
  \begin{align*}
   \Omega_{q*}(\mathcal{G},\epsilon) \cap \Omega_q(\mathcal{G},\epsilon)
   \subset \Ospp(\epsilon)\;.
\end{align*} 
\end{proposition}
In words, this proposition states that if the first $J_0$ of the
$G^{(j)}$ are small and all of the quantities
$\N{X^{(j)}}_{1},\N{Y_i^{(j)}}_{1}$ are not to big then it is unlikely
that the remaining $G^{(j)}$ are big.

\bpf[Proof of Proposition \ref{l:ladderEstimateLocalized}] The idea of
the proof is to iterate Lemma \ref{l:oneStep} below. Begin by setting
$g_*=\sup_{j} \sum_{k=1}^{K(j)} |g_k|_\infty$.  Now define
$\epsilon_0=\tilde \epsilon_0=\epsilon$ and for $j > 1$ define
$\epsilon_{j}=\tilde\epsilon_{j-1}^\frac1{151}$, $\tilde
\epsilon_j=\epsilon_j^\frac{151}{152}$. With these choices,
$\epsilon_j < \tilde \epsilon_j < \epsilon_{j+1} < \tilde
\epsilon_{j+1}$. We also choose $q_0$ sufficiently small so that
$\tilde \epsilon_{J+1} < \epsilon^q$ for $\epsilon \in [0,1)$ and $q
\in(0,q_0]$. Define the following subsets of $\Omega$
\begin{align*}
  A(j)=&  \Omega_{q*}(\mathcal{G},\epsilon) \cap \Omega_q(\mathcal{G},\epsilon) \cap \left\{ \sup_{1\leq i\leq
      N}\big(\|X^{(l)}\|_{\infty},\|Y_i^{(l)}\|_{\infty}\big)
    \leq  \epsilon_{l+1} \mbox{ for $l \leq j$} \right\}  \\
  A^+(j)=& A(j) \cap \left\{ \|G^{(l)}\|_\infty \leq \tilde
    \epsilon_{l} \leq \epsilon_{l+1}\mbox{ for $l \leq j+1$} \right\}  \\
  B(j)=&\left\{ \sup_{ 1\leq i\leq
      N}\big(\|X^{(j)}\|_{\infty},\|Y_i^{(j)}\|_{\infty}\big)>
    \tilde\epsilon_j^\frac1{151}= \epsilon_{j+1} \right\} \\
  C(j)=& \Omega_{q*}(\mathcal{G},\epsilon) \cap \left\{
    \|G^{(j)}\|_{\infty} \leq \tilde\epsilon_j \ ;\sup_{ 1\leq i\leq
      N}\big(\|X^{(j)}\|_{\infty},\|Y_i^{(j)}\|_{\infty}\big) >
    \tilde\epsilon_j^\frac1{151}= \epsilon_{j+1} \right\}.
\end{align*}

First notice that since $\epsilon < \epsilon^q$ for $\epsilon \in
(0,1]$ and $q\in (0,q_0]$, the event
\begin{align*}
    \biggl\{ \sup_{1\leq j\leq J_0}
   \|G^{(j)}\|_{\infty} < \epsilon \ ;\sup_{1\leq j\leq J_0}\|G^{(j)}\|_\infty>\epsilon^q \biggr\}
\end{align*}
is empty. Hence \begin{align*} 
    \Omega_q(\mathcal{G},\epsilon)&=\Biggl\{ \sup_{1\leq j\leq J_0}
   \|G^{(j)}\|_{\infty} < \epsilon \ ;\sup_{\substack{1\leq j\leq J\\
       1\leq i\leq N\\ J_0 < \ell\leq J}}\big(\|X^{(j)}\|_{\infty},
   \|Y_i^{(j)}\|_{\infty}, \|G^{(\ell)}\|_\infty\big)>\epsilon^q \Biggr\}
   \; .
  \end{align*} 
Next notice that for any $j>J_0$, because $G^{(j)}(t)=\sum_{k=1}^K
g_k(t)f_k(t)$, $\|G^{(j)}\|_\infty \leq g_* \sup_k \|f_k\|_\infty$
where the $f_k$ are from earlier in the ladder in the sense of
\eqref{eq:earlyInLadder}. Hence for any $j$, we have that if
\begin{align*}
  \sup_{1\leq i\leq N}\big(\|X^{(l)}\|_{\infty},
  \|Y_i^{(l)}\|_{\infty}, \|G^{(l)}\|_{\infty}\big) \leq
  \epsilon_{l+1} \mbox{ for $l \leq j-1$}
\end{align*}
then $\|G^{(j)}\|_\infty \leq g_* \epsilon_{j} <
\epsilon_{j}^\frac{151}{152}= \tilde \epsilon_{j}$ for $\epsilon$
sufficiently small. Restricting to $\epsilon$ small enough so that $
g_* \epsilon_{j} < \epsilon_{j}^\frac{151}{152}$ for all $1\leq j \leq
J$ implies that $A^+(j) =A(j)$ for all $1\leq j \leq J$.

Next observe that $A(j) \cap B(j+1)^c = A(j+1)$ and hence $A^+(j) \cap
B(j+1)^c = A^+(j+1)$ for $\epsilon$ sufficiently small. Iterating this
observation, with the convention $A^+(0)=
\Omega_{q*}(\mathcal{G},\epsilon) \cap
\Omega_q(\mathcal{G},\epsilon)$, we obtain
\begin{align*}
  A^+(0)&=  [A^+(0)\cap B(1)^c] \cup [A^+(0)\cap B(1)] \\&= A^+(1)
  \cup  [A^+(0)\cap B(1)]\\
  &= [A^+(1)\cap B(2)^c] \cup [A^+(1)\cap B(2)] \cup [A^+(0)\cap
  B(1)]\\
  &=A^+(2) \cup  [A^+(1)\cap B(2)] \cup  [A^+(0)\cap B(1)]\\
  &= A^+(J) \cup \bigcup_{j=0}^{J-1}\big[ A^+(j)\cap B(j+1) \big] \ .
\end{align*}
Since  $q$ was picked sufficiently small so that $\epsilon_{J+1} < \epsilon^q$,
we observe that $A^+(J)$ is empty since on $A^+(J)$ 
\begin{align*}
  \sup_{\substack{ 1\leq j \leq J\\1\leq i\leq N \\J_0< l\leq
      J}}\big(\|X^{(j)}\|_{\infty},\|Y_i^{(j)}\|_{\infty},\|G^{(l)}\|_{\infty}\big)
  \leq \epsilon_{J+1} <\epsilon^q< \sup_{\substack{ 1\leq j \leq
      J\\1\leq i\leq N\\ J_0< l\leq
      J}}\big(\|X^{(j)}\|_{\infty},\|Y_i^{(j)}\|_{\infty},\|G^{(l)}\|_{\infty}\big)
\end{align*}
which cannot be satisfied.  

Recall that $\Ospp(\epsilon)=\bigcup_{j=1}^J \Os(\tilde \epsilon_j)$.
Let $ \OII(H^{(j)},\tilde \epsilon_j)$ be the set defined bellow in
Lemma \ref{l:oneStep}.  For all $q$ sufficiently small,
$\Omega_{q*}(\mathcal{G},\epsilon) \subset \OII(H^{(j)},\tilde
\epsilon_j)$ for $j=1,\cdots,J+1$.  Decrease $q_0$ so this holds.
With this choice Lemma \ref{l:oneStep} implies that $C(j) \subset
\Os(\tilde \epsilon_j)$.  Since clearly $A^+(j) \cap B(j+1)
\subset C(j+1)$, combining all of these observations produces
\begin{align*}
  \Omega_{q*}(\mathcal{G},\epsilon) \cap
  \Omega_q(\mathcal{G},\epsilon)= A^+(0)\subset \bigcup_{j=1}^J C(j)
  \subset \bigcup_{j=1}^J \Os(\tilde \epsilon_j) =\Ospp(\epsilon)\;.
\end{align*}
\epf

Lastly we give a version of the preceding proposition which begins with $L^p$
estimates in time on the $\{G^{(j)} : 1\leq j \leq J_0\}$ rather than
$L^\infty$ estimates.
\begin{corollary} \label{c:LpLadder} Fix $T>0$. For any $\ell >0$,  define 
  \begin{equation*}
   \Omega_{q,\ell}(\mathcal{G},\epsilon)=\Biggl\{ \sup_{1\leq j\leq
        J_0} \int_0^T \big|G^{(j)}(s)\big|^\ell ds < \epsilon \ 
      ;\sup_{\substack{1\leq j\leq J\\ 1\leq i\leq N}}\big(\|X^{(j)}\|_{\infty}, \|Y_i^{(j)}\|_{\infty},
      \|G^{(j)}\|_\infty\big) > \epsilon^q \Biggr\}\;.
  \end{equation*}
  There exist positive constants $q=q(J,\ell)$ and
  $\epsilon_0=\epsilon_0(J,T,N,\ell)$ so that for all $\epsilon
  \in(0,\epsilon_0]$
    \begin{equation*}
      \Omega_{q*}(\mathcal{G},\epsilon) \cap
      \Omega_{q,\ell}(\mathcal{G},\epsilon) \subset 
      \Ospp(\epsilon)\;. 
\end{equation*}
\end{corollary}
\bpf  We begin by translating the bound
\begin{equation*}
  \sup_{1\leq j\leq J_0} \int_0^T \big|G^{(j)}(s)\big|^\ell ds <
\epsilon
\end{equation*}
into a bound of the form $\sup_{1\leq j\leq J_0}
\|G^{(j)}\|_\infty \leq \epsilon^\beta$ for some $\beta \in (0,1)$.  

Notice that 
\begin{align*}
\left|G^{(j)}(s)-  G^{(j)}(r)\right| = \left| \int_r^s H^{(j)}(t)
  dt\right| \leq |s-r| \;\|H^{(j)}\|_{\infty}\;.
\end{align*}
Without loss of generality we assume $q < \frac1{150}$. Hence on
$\Omega_{q*}(\mathcal{G},\epsilon)\cap \left\{ \sup_i \|W_i\|_\infty >
  \epsilon^{-\frac{1}{150}} \right\}^c$
\begin{align*}
  \|H^{(j)}\|_{\infty} \leq \|X^{(j)}\|_{\infty} + \sum_{i=1}^N
  \|Y_i^{(j)}\|_{\infty} \|W_i^{(j)}\|_{\infty} \leq
  \epsilon^{-\frac{1}{150}} + N \epsilon^{-\frac{1}{75}} \leq (N+1)
  \epsilon^{-\frac{1}{75}}
\end{align*}
In other words, $\hol_1(G) \leq(N+1) \epsilon^{-\frac1{75}}$ on $
\Omega_{q*}(\mathcal{G},\epsilon)\cap \left\{ \sup_i \|W_i\|_\infty >
  \epsilon^{-\frac1{150}} \right\}^c$. Then Lemma \ref{l:PToInfinty} below\,
implies $\sup_{j\leq J_0}\|G^{(j)}\|_{\infty} \leq
(2+N)\epsilon^{\beta_0} < \epsilon^{\beta_1}$ for $\epsilon$
sufficiently small where $\beta_0=\frac{74}{75}\frac{1}{1+\ell}$ and
$\beta_1=\frac{1}{2}\frac{1}{1+\ell}$. Now notice that
$\Omega_*(\mathcal{G},\epsilon^{\beta_1})\subset
\Omega_*(\mathcal{G},\epsilon)$ because $\beta_1 < 1$. Hence the
result follows from Proposition \ref{l:ladderEstimateLocalized} and
the fact that $\left\{\sup_i \|W_i\|_\infty >
  \epsilon^{-\frac1{150}} \right\} \subset \Ospp(\epsilon)$.  \epf
\subsection{The Basic Estimates}
Let
\begin{align*}
  G(t)=G_0+ \int_0^t H(s) ds
\end{align*}
Now let $H(s)$ be any stochastic process of the form
\begin{equation}\label{eq:H}
  H(s)=X(s)+\sum_{i=1}^N Y_i(s)W_i(s)\eqdef X(s)-Z(s)
\end{equation}
where $X(s)$, and $Y_1(s),\ldots,Y_N(s)$ are Lipschitz continuous stochastic
processes and $\{W_1(s),$ $\ldots,W_N(s)\}$ are mutually independent standard Wiener 
processes with
$W_i(0)=0$, $1\le i\le N$. 

Next given $\epsilon>0$, define the following subsets of the  probability space:
\begin{align}\label{eq:boundsLipSize}
  \OI(H,\epsilon) &= \left\{ \sup_{1\leq i\leq N} \big(
    \|Y_i\|_\infty \big) \leq \epsilon^{-\frac1{28}} \ ; \ 
    \sup_{1\leq i\leq N} \big( \hol_1(Y_i), \hol_1(X) \big) \leq
    \epsilon^{-\frac1{28}} \right\} \;,\\
    \OII(H,\epsilon)&=\left\{ \sup_{1\leq i\leq N}\big(\N{X}_{1},\N{Y_i}_{1}\big)\leq
      \epsilon^{-\frac1{150}} \right\}\;.  \notag  
\end{align}

\begin{lemma}\label{l:GtoH} Let $\epsilon >0$.
  Assume $\hol_\alpha(H) \leq c\epsilon^{-\gamma}$ for some fixed
  $\alpha > \gamma >0$. Then $\|G\|_{\infty} \leq \epsilon$ implies
  $\|H\|_{\infty} \leq (2+c)\epsilon^\frac{\alpha-\gamma}{1+\alpha}$.
\end{lemma}
\bpf For any $s \in [0,t]$, let $r_1 \leq s \leq r_2$ such that
$r_2-r_1=\epsilon^{\frac{1+\gamma}{1+\alpha}}$.  Notice that by the
assumption $H(r) \geq H(s) - c\epsilon^{-\gamma} |r_2-r_1|^\alpha$ for
any $r \in [r_1,r_2]$. Hence we have
\begin{align*}
  2\epsilon \geq G(r_2)-G(r_1) = \int_{r_1}^{r_2} H(r) dr \geq
  |r_2-r_1|( H(s) -  c\epsilon^{-\gamma}
|r_2-r_1|^\alpha) \ .
\end{align*}
Rearranging this gives, $H(s) \leq \frac{2\epsilon}{|r_2-r_1|} +
c\epsilon^{-\gamma} |r_2-r_1|^\alpha= (2+c)
\epsilon^\frac{\alpha-\gamma}{1+\alpha}$. The same argument from above
gives a complementary lower bound and completes the result.  \epf

Next, we have the following result.
\begin{lemma} \label{l:oneStep}
  There exists a $\epsilon_0=\epsilon_0(N)$ so that for every
  $\epsilon \in (0,\epsilon_0]$ and stochastic process $G(t)$, of the
  form given above, one has
  \begin{align*}
\OII(H,\epsilon)   \cap \left\{ \|G\|_{\infty} < \epsilon \ ; \
      \sup_{1\leq i\leq N}
      \big(\|X\|_{\infty},\|Y_i\|_{\infty} \big)
      > \epsilon^\frac1{151}\right\} \subset
    \Os(\epsilon)\;.
  \end{align*}
\end{lemma}
\bpf The result will follow from Lemma \ref{l:likelyNotSmall} of the
next subsection after some ground work is laid.  As before, set $\Delta=
\epsilon^\frac{14}{75}$ and recall that by definition $\Os(\epsilon)=
\Omega_a(\Delta)\cup\Omega_b(\Delta)\cup\Omega_c(\Delta)$. Then notice
that $\hol_\frac14(f) \leq \hol_1(f)$, $\hol_\frac14(f+g) \leq
\hol_\frac14(f)+ \hol_\frac14(g)$, and $\hol_\frac14(fg)\leq
\hol_\frac14(g)\|f\|_{\infty} +\hol_\frac14(f)\|g\|_{\infty}$ for all functions
$f$ and $g$. Hence on $\Oc(\Delta )^c \cap \OII(H,\epsilon)$
\begin{align*}
  \hol_\frac14(H) &\leq \hol_1(X) + \sum_i \hol_1(Y_i)\|W_i\|_{\infty}
  +
  \hol_\frac14(W_i)\|Y_i\|_{\infty} \\
  & \leq \epsilon^{-\frac1{150}} +2 N \epsilon^{-\frac1{75}} \leq
  (1+2N)\epsilon^{-\frac1{75}} \ .
\end{align*}
Hence by Lemma \ref{l:GtoH} one has that on $\Oc(\Delta )^c \cap
\OII(H,\epsilon)$,  $\|G\|_{\infty} < \epsilon$ implies $\|H\|_{\infty} \leq
(1+2N)\epsilon^\frac{71}{375} <\Delta=\epsilon^\frac{14}{75}$ for all
$\epsilon$ sufficiently small.  Next observe that because $\Delta=
\epsilon^\frac{14}{75}$, $\OII(H,\epsilon)\subset \OI(H,\Delta)$ where
$\OI(H,\Delta)$ is the set which was defined in equation \eqref{eq:boundsLipSize}. In light of this, Lemma \ref{l:likelyNotSmall}
implies that
\begin{equation}\label{eq:inclusionPrevious}
  \OII(H,\epsilon)\cap \left\{ \|G\|_{\infty} < \epsilon \ ; \ 
    \sup_i\|Y_i\|_{\infty} > \epsilon^\frac1{75} \right\} \subset
  \Os(\epsilon)\;.
\end{equation}
Now on 
\begin{align*}
 \Oc(\Delta )^c \cap \OII(H,\epsilon)\cap \left\{ \|G\|_{\infty} < \epsilon \ ; \ 
    \sup_i\|Y_i\|_{\infty} \leq \epsilon^\frac1{75} \right\}
\end{align*} one has that
$\|X\|_{\infty} \leq \|H\|_{\infty} + \sum
\|Y_i\|_{\infty}\|W_i\|_{\infty} \leq \epsilon^\frac{14}{75} + N
\epsilon^\frac1{75}\epsilon^{-\frac1{150}}$ which is less than
$\epsilon^\frac1{151}$ for $\epsilon$ sufficiently small. Hence 
\begin{align*}
 \Oc(\Delta )^c \cap \OII(H,\epsilon)\cap \left\{ \|G\|_{\infty} < \epsilon \ ; \ 
    \sup_i\|Y_i\|_{\infty} \leq \epsilon^\frac1{75} \right\}
    \cap\left\{  \|X\|_{\infty} > \epsilon^\frac1{151}   \right\}
\end{align*} 
is empty. Combining this observation with \eqref{eq:inclusionPrevious}
implies
\begin{align*}
  \OII(H,\epsilon)\cap \left\{ \|G\|_{\infty} < \epsilon \right\} \cap
    \left[  \Bigl\{ \sup_i\|Y_i\|_{\infty} > \epsilon^\frac1{75}\Bigr\}\cup\left\{
    \|X\|_{\infty} > \epsilon^\frac1{151} \right\}\right] \subset
  \Os(\epsilon)\;.
\end{align*}
Since $\{ \sup_i\|Y_i\|_{\infty} > \epsilon^\frac1{151}\} \subset \{
 \sup_i\|Y_i\|_{\infty} > \epsilon^\frac1{75}\} $, the proof
is complete.\epf

\begin{lemma} \label{l:PToInfinty}For any $\epsilon>0$ and $\ell>0$,
  $\int_0^T |G(s)|^\ell ds < \epsilon$ and $\hol_\alpha(G) <c
  \epsilon^{-\gamma}$ implies $\|G\|_{\infty} <(1+c)
  \epsilon^\frac{\alpha-\gamma}{1+\ell\alpha}$
\end{lemma}
\bpf By Chebyshev's inequality for any $\beta >0$, we have $\lambda\{
x: |G(x)| > \epsilon^\beta\} \leq \epsilon^{1-\ell\beta}$ where $\lambda$
is Lebesgue measure. Hence $\|G\|_{\infty} \leq \epsilon^\beta +
c\epsilon^{(1-\ell\beta)\alpha} \epsilon^{-\gamma}$.
Setting $\beta=\frac{\alpha-\gamma}{1+\ell\alpha}$ proves the result.
\epf

\subsection{The Main Technical Estimate: }
Let $H(s)$ be as in \eqref{eq:H} from the preceding section.
Define the following piecewise constant approximation of the $Z$ from the
definition of $H$:
\begin{equation}\label{eq:Z*}
  Z^*(s)=-\sum_{i=1}^N Y_i^*(s)W_i(s)
\end{equation}
where $Y_i^*(s)= \sum_{k=1}^m \ONE_{I_k}(s) \sup_{s\in I_k}Y_i(s)$
and $\ONE_{I_k}$ is the indicator function of the set $I_k=[t_{k-1},t_k)$.

For any $k$, and process $\zeta(s)$ defined on $[t_{k-1},t_k]$, we
define the $\delta$-scale quadratic variation on $[t_{k-1},t_k]$ by 
\begin{align*}
  Q_k(\zeta)=\sum_{\ell=1}^{M(k)} \bigl[ \delta_\ell^k \zeta\bigr]^2 \ .
\end{align*}

\begin{lemma} \label{l:mainDeterministicSmall} On $\OI(H,\Delta) \cap
  \Omega_c(\Delta)^c$ for any $k=1,\dots,m$, we have the following estimates:
  $Q_k(X) \leq 2\Delta^{\frac{109}{42}}$, $Q_k(Z^*)\leq
  2Q_k(Z)+4N^2\Delta^{\frac{25}{21}}$. Moreover, on $\OI(H,\Delta) \cap
  \Omega_c(\Delta)^c\cap\left\{\sup_{s \in I_k} |H(s)| < \Delta\right\}$,
  $Q_k(Z)\leq 2 Q_k(X) + 8\Delta^{\frac43}+8\Delta^2$ and 
  $Q_k(Z^*) \leq (40+4N^2) \Delta^{\frac{25}{21}}$.
\end{lemma}
\bpf

For brevity, we suppress the $k$ dependence of $M(k)$ and
$s_\ell(k)$. The first inequality follows from
\begin{align*}
  Q_k(X)&=\sum_{\ell=1}^M  [X(s_\ell) -X(s_{\ell-1})]^2 \\
  &\leq \sum_{\ell=1}^M [\delta \Delta^{-\frac{1}{28}}]^2 \leq
  M\delta^2 \Delta^{-\frac{1}{14}} =
  (\Delta^{-\frac23}+1)\Delta^{\frac{10}{3}}\Delta^{-\frac{1}{14}}=
  \Delta^{\frac{109}{42}}+ \Delta^{\frac{139}{42}}< 2 \Delta^{\frac{109}{42}} \; .
\end{align*}
To see the second implication, first notice that
\begin{align*}
  Q_k(Z^*) \leq& \sum_{\ell=1}^M \left(
    2[|Z^*(s_{\ell-1})-Z(s_{\ell-1})| + |Z^*(s_{\ell})-Z(s_{\ell})|
    ]^2 +2[Z(s_{\ell-1})-
    Z(s_{\ell})]^2 \right)\\
  \leq & 2Q_k(Z) + 4 \sum_{\ell=0}^{M-1}
  [Z^*(s_{\ell})-Z(s_{\ell})]^2+4 \sum_{\ell=1}^{M}
  [Z^*(s_{\ell})-Z(s_{\ell})]^2.
\end{align*}
Next note that
\begin{align*}
  \sum_{\ell=1}^M [Z^*(s_{\ell})-Z(s_{\ell})]^2 &\leq \sum_{\ell=1}^M
  \left[\sum_{i=1}^N[Y_i(s_\ell)-Y_i^*]W_i(s_\ell)\right]^2 \\ 
  &\leq \Delta^{-\frac17}\frac{N^2}{\delta} \sum_{\ell=1}^M (\ell
  \delta)^2 \delta \leq \frac{N^2}{3\delta}\Delta^{\frac{20}{7}} =
  \frac{N^2}{3}\Delta^{\frac{25}{21}}
\end{align*}
and similarly $\sum_{\ell=0}^{M-1} [Z^*(s_{\ell})-Z(s_{\ell})]^2 \leq
\frac{N^2}{3}\Delta^{\frac{25}{21}} $. Combining this estimate with the
previous gives the second result. 

The third result follows from $H(s)=X(s)-Z(s)$ and
\begin{align*}
  Q_k(Z) &= \sum_{\ell=1}^M [Z(s_{\ell})-Z(s_{\ell-1})]^2\\
  &\leq \sum_{\ell=1}^M 2[|Z(s_{\ell})-X(s_{\ell})| +
  |Z(s_{\ell-1})-X(s_{\ell-1})| ]^2 +2[X(s_{\ell})-X(s_{\ell-1})]^2\\
  &\leq 8\Delta^{2} M + 2Q_k(X)\leq 8\Delta^2(\frac{\Delta}{\delta}+1)
  + 2Q_k(X)\; .
\end{align*}
Lastly, combining the three previous estimates produces the final estimate.
\epf

To aid in the analysis of $Y^*$, consider a general process of the form
\begin{align*}
  \zeta(s)=\sum_{i=1}^N a_i(s) W_i(s)
\end{align*}
where the $W_i$ are independent standard Wiener processes and the
$a_i(s)$ are constant on the intervals $I_k=[(k-1)\Delta,k\Delta)$
for each $k=1,\dots,m$. As before, for $k=1,\dots m$, we define
\begin{align*}
  Q_k(\zeta)=\sum_{\ell=1}^{M(k)} \Biggl( \sum_{i=1}^N
  a_i(s_\ell(k)) \delta_\ell^kW_i \Biggr)^2 \;,
\end{align*}
where $s_\ell(k)$ and $\delta_\ell^k$ are as defined at the start of
Section \ref{sec:ladderDef}. Notice that if we define
\begin{xalignat*}{3}
  U&=\sum_{i=1}^N a_i^2\sum_{\ell=1}^{M(k)} (\hat\delta_\ell^kW_i)^2 &
  V&=\sum_{\substack{(i,j) \in \{1,\dots,N\}^2 \\ i\neq j}} a_i a_j
  \sum_{\ell=1}^{M(k)}(\hat\delta_\ell^kW_i)(\hat\delta_\ell^kW_j) ,
\end{xalignat*}
where $a_i=a_i\big((k+1)\Delta\big)$ and $\hat \delta_\ell^k f=(\delta_\ell^k
f)/\sqrt{\delta_\ell^k}$ then $Q_k(\zeta)=\frac{\Delta}{M(k)}(U+V)$.

\begin{lemma}\label{l:mainProbSmall}For $\sigma > \frac{8}{7}$, $\Delta \in(0,6^{\frac{-7}{7\sigma-8}})$ and
  $k=1,\dots,m$,
  \begin{align*}
    \left\{ Q_k(\zeta) < \Delta^{\sigma} \ ;  \Delta^\frac{1}{14} < \sup_i |a_i| \leq
      \Delta^{-\frac1{28}}\right\} \subset \Omega_a(\Delta) \cup
    \Omega_b(\Delta)
    \end{align*}
\end{lemma}
\bpf First notice that because $\frac{8}{7} < \sigma$ and $\Delta <
6^\frac{-7}{7\sigma-8}$, $\Delta^{\sigma}- \frac12\Delta^\frac87 < -
\frac13\Delta^\frac87$ so
\begin{align*}
  &\left\{ Q_k(\zeta) < \Delta^{\sigma} \ ;  \Delta^\frac{1}{14} < \sup_i |a_i| \leq
      \Delta^{-\frac1{28}}\right\} \\&\quad \subset \left\{ \frac{\Delta}{M(k)}U <
    \frac12\Delta^\frac87 \ ;\ \sup_i |a_i| >
    \Delta^\frac{1}{14}\right\} \cup \left\{ \frac{\Delta}{M(k)}V <
    \Delta^{\sigma} - \frac12\Delta^\frac87
    \ ;\ \sup_i |a_i| < \Delta^{-\frac1{28}}\right\} \\
  &\quad\subset \left\{ \frac{\Delta}{M(k)}U < \frac12\Delta^{\frac87}
    \ ;\ \sup_i |a_i| > \Delta^\frac1{14}\right\} \cup \left\{
    \frac{\Delta}{M(k)}|V| > \frac13\Delta^{\frac87}
    \ ;\  \sup_i |a_i| < \Delta^{-\frac1{28}}\right\} \\
\end{align*}
Now 
\begin{align*}
  \Bigl\{ U \leq \frac12\Delta^\frac17 M(k) \ ; \ \sup_i |a_i| >
  \Delta^\frac{1}{14}\Bigr\} & \subset \Bigl\{ \inf_i
  \sum_{\ell=1}^{M(k)} (\hat \delta_\ell^k W_i)^2 \leq \frac12 M(k)\Bigr\} \subset \Omega_a(\Delta)
\end{align*}
and
\begin{multline*}
  \left\{ |V| > \frac13\Delta^\frac17 M(k) \ ;\ \sup_i |a_i| <
    \Delta^{-\frac{1}{28}} \right\} \\\subset \Biggl\{
  \sup_{\substack{(i,j)\\ i\neq j}} |a_i| |a_j| \left|
    \sum_{\ell=1}^{M(k)}(\hat\delta_\ell^kW_i)(\hat\delta_\ell^kW_j)\right|
  > \frac{\Delta^\frac17}{3N^2} M(k) \ ;\ \sup_i |a_i| <
  \Delta^{-\frac1{28}} \Biggr\} \\\subset \Biggl\{
  \sup_{\substack{(i,j)\\ i\neq j}} \left|
    \sum_{\ell=1}^{M(k)}(\hat\delta_\ell^kW_i)(\hat\delta_\ell^kW_j)\right|
  > \frac{\Delta^{\frac{3}{14}}}{3N^2} M(k) \Biggr\} \subset \Omega_b(\Delta)
\end{multline*}
\epf

The following result is the main result of this section.
\begin{lemma} \label{l:likelyNotSmall}  
 For all $\Delta \in (0,(40+4N^2)^{-42}]$ and all stochastic processes
  $H(s)$ of the form \eqref{eq:H}
\begin{equation*}
 \OI(H,\Delta) \cap \Big\{\|H\|_{\infty} < \Delta \ ; \
   \sup_{i}\|Y_i\|_{\infty} > \Delta^\frac1{14} \Big\} 
\subset \Omega_a(\Delta)\cup\Omega_b(\Delta)\cup\Omega_c(\Delta)\;.
\end{equation*}
\end{lemma}
We will prove this result by showing that on $\OI(H,\Delta)$ as
$\Delta \rightarrow 0$, if
\begin{equation*}
\sup_{s \in [0,T]} |H(s)| < \Delta \mbox{ and }\sup_{i}\sup_{s
  \in [0,T]} |Y_i(s)| > \Delta^\frac{1}{14}  
\end{equation*}
then the approximate quadratic variation of the Wiener processes at
the scale $\delta$ is abnormally small or $\sup_i \|W_i\|_\infty >
\Delta^{-\frac1{28}}$.

\bpf[Proof of Lemma \ref{l:likelyNotSmall}] From the last estimate in
Lemma \ref{l:mainDeterministicSmall}, we have that on $\OI(H,\Delta)
\cap \Omega_c(\Delta)^c\cap \{\sup_{s \in I_k} |H(s)| < \Delta \}$,
$Q_k(Z^*) \leq (40+4N^2) \Delta^\frac{25}{21} <\Delta^\frac{49}{42}$
for $\Delta \in (0,(40+4N^2)^{-42}]$. Here $Z^*$ is the approximation
defined in \eqref{eq:Z*}.  Now Lemma \ref{l:mainProbSmall} with
$\sigma=\frac{49}{42}( > \frac87)$ implies that
\begin{equation}\label{eq:mainLocalEstimate}
  \Omega_c(\Delta)^c\cap \OI(H,\Delta)\cap \left\{ \sup_{s \in I_k}
  |H(s)| < \Delta \ ; \ 
  \sup_{i}\sup_{s \in  I_k} |Y_i(s)| > \Delta^\frac{1}{14} \right\}
  \subset  \Omega_a(\Delta) \cup \Omega_b(\Delta) 
\end{equation}
for all $k=1,\dots$ and $\Delta \in (0,(40+4N^2)^{-42}]$.

Continuing, we have that
\begin{align*}
  \Omega_c(\Delta)^c&\cap \OI(H,\Delta)\cap \left\{ \|H\|_{\infty} <
    \Delta ; 
    \sup_{i} \|Y_i\|_{\infty} > \Delta^\frac{1}{14}\right\}\\
  &\subset \bigcup_{k=1}^m \left\{ \sup_{s \in I_k} |H(s)| < \Delta \ 
    ;\ \sup_{i}\sup_{s \in I_k} |Y_i(s)| > \Delta^\frac{1}{14}
  \right\}
  \cap   \OI(H,\Delta) \cap \Omega_c(\Delta)^c\\
  &\subset  \Omega_a(\Delta) \cup \Omega_b(\Delta) 
\end{align*}
\epf
\subsection{Estimates on the Size of $\Omega_a$, $\Omega_b$, and
  $\Omega_c$}\label{sec:SmallABC}
Since the events described by  $\Omega_a$ and $\Omega_b$ are simply
statements about collections of independent standard normal random
variables, the following two estimates will give us the needed control.
\begin{lemma} \label{l:NormalProb} For $c \in  (0,1)$ and $M >
  \frac2{1-c}$, setting $\gamma=c-1-\ln(c)>0$
  \begin{align*}
    \P \Bigl( \sum_{\ell=1}^M \eta_{\ell}^2 \leq c M \Bigr)
    &\leq \frac1{\sqrt{\pi M}}\exp\Big( - \frac12\gamma M\Big)\\
    \P \Bigl( \Bigl|\sum_{\ell=1}^M \eta_\ell \tilde\eta_\ell
    \Bigr|\geq cM\Bigr) & \leq 2\P \Bigl( \sum_{\ell=1}^M \eta_\ell
    \tilde\eta_\ell \geq cM\Bigr) \leq
    2\exp\Bigl(-\frac{c^2}{4}M\Bigr)
  \end{align*}
where $\{ \eta_{\ell},\tilde\eta_\ell\}$ are a collection of $2M$
mutually independent standard  $N(0,1)$ random
variables.
\end{lemma}
\bpf
Notice that  $\sum_{\ell=1}^M \eta_{\ell}^2$ is distributed as
a $\chi^2$ random variable with $M$ degrees of freedom. Hence we have
\begin{align*}
 \P \Bigl(\sum_{\ell=1}^M \eta_{\ell}^2 \leq c M\Bigr) &=
 \frac{2^{-\frac{M}{2}}}{\Gamma(\frac{M}{2})} \int_0^{c M}
 x^{\frac{M}{2}-1} e^{-\frac{x}{2}} dx
\end{align*}
Since $c <1$ and $M>\frac2{1-c}$, the integrand is bounded by
$(cM)^{\frac{M}2-1} \exp(-c \frac{M}{2})$. Combining
this with $\Gamma(\frac{M}{2}) \geq \sqrt{\pi M}
(\frac{M}{2e})^\frac{M}2$ implies that  
\begin{align*}
 \P \Bigl(\sum_{\ell=1}^M \eta_{\ell}^2 \leq c M\Bigr) &\leq
 \frac{M^{-\frac{M}{2}}e^{\frac{M}{2} } }{\sqrt{\pi M}}
 (cM)^{\frac{M}2} \exp(-c\frac{M}{2})\\
 &\leq \frac1{\sqrt{\pi M}} \exp( \frac{M}{2} [ -c  + 1
 +\ln(c)]) \ .
\end{align*}
Noticing that $-c + 1 +\ln(c)< 0$ for $c \in (0,1)$ finishes the
proof of the first statement. 

For the second, note that for $\lambda \in (-1, 1)$, $\E \exp(\lambda
\eta_\ell \tilde \eta_\ell)= (1-\lambda^2)^{-\frac12}$. Hence for $\lambda
\in[0,\frac12]$, 
\begin{align*}
  \P \Bigl( \sum_{\ell=1}^M
      \eta_\ell \tilde\eta_\ell \geq cM\Bigr) &\leq \exp\left(
      -\lambda c M + \frac12|\ln(1-\lambda^2)|M\right)\\ &\leq \exp\left(
      -\lambda c M + \lambda^2 M\right)
\end{align*}
Taking $\lambda =\frac{c}{2}$ gives the result.
 \epf
 \begin{corollary} \label{c:ab}For $\Delta \leq \frac12 \wedge T$,
   \begin{align*}
     \P\Bigl( \Omega_a(\Delta) \Bigr) &\leq \frac{2TN}{\sqrt{\pi}}
     \Delta^{-\frac23} \exp\Bigl(-\frac1{20} \Delta^{-\frac23}\Bigr)\\
     \P\Bigl(\Omega_b(\Delta) \Bigr) & \leq 6 N^2 T \Delta^{-1}
     \exp\Bigl(- \frac1{3N^2} \Delta^{-\frac{19}{42}}\Bigr)\;. 
   \end{align*}
   In particular, if $\gamma=\min(\frac1{3N^2},\frac1{20})$ then
   \begin{align*}
     \P\Bigl(  \Omega_a(\Delta) \cup  \Omega_b(\Delta) \Bigr) \leq 8
     T N^2 \frac{\exp(-\gamma \Delta^{-\frac25})}{\Delta}
   \end{align*}
 \end{corollary}
\bpf
First observe that the $\{\hat \delta_\ell^kW_i\}$ are independent $N(0,1)$
random variables. Since
\begin{align*}
  \Omega_a\subset \bigcup_{k=0}^m \bigcup_{i=1}^N \left\{
  \sum_{\ell=1}^{M(k)}  (\hat\delta_\ell^kW_i)^2 \leq \frac{M(k)}2
  \right\}\; , 
\end{align*} 
$m \leq \frac{T}{\Delta}+1$ and $\Delta^{-\frac23}-1\leq M \leq
\Delta^{-\frac23}$ the first result follows from Lemma
\ref{l:NormalProb} and bounding the previous expression by the sum of
the probability of the sets on the right handside. 

Proceeding in a fashion similar to the first estimate, the second bound
follows from
\begin{align*}
  \Omega_b\subset \bigcup_{k=0}^m \bigcup_{\substack{(i,j) \in
      \{1,\dots,N\}^2\\ i\neq j}} \left\{ \left|\sum_{\ell=1}^{M(k)}
      (\hat\delta_\ell^kW_i) (\hat\delta_\ell^kW_j) \right|\geq
    \frac{\Delta^\frac{3}{14}M(k)}{3N^2} \right\}\; ,
\end{align*} 
Combining the first two estimates gives the last quoted result.
\epf
\begin{lemma}\label{l:c} For any $p\geq 1$ there exists a $c=c(p,N,T)$
  so that $\P\big(\Omega_c(\epsilon) \big) \leq c \epsilon^p$.   
\end{lemma}
\bpf It is enough to show that $\E \|W_i\|_\infty^\gamma <
c(T,\gamma)$ and $\E \big[\hol_\frac14(W_i)^\gamma\big] <
c(T,\gamma)$ for any $\gamma \geq 1$. The first follows from the
Doob's inequalities for the continuous martingale
$W_i(s)$. The finiteness of the moments of the modules continuity of
Wiener processes is given by Theorem 2.1 (p.  26) and the observation
at the top of p.  28 both in \cite{b:RevuzYor94}. Together these imply
that $\E \big[\hol_\alpha(W_i)^\gamma\big] < \infty$ for all
$\gamma >0$ as long as $\alpha \in (0,\frac12)$. Since
$\alpha=\frac14$ in our setting, the proof is complete.  \epf

\section{Strict Positivity of the Density}
\label{sec:positive}

We will now give conditions under which for any $t>0$ and some
orthogonal projection $\Pi$ of $\L^2$ onto a finite dimensional
subspace $S$, the density $p(t,x)$ of the law of $\Pi \sol(t)$ with
respect to Lebesgue measure on $S$ satisfies
\begin{align}
  \label{eq:pos}
  p(t,x) > 0\;, \text{for all $x \in S$.}
\end{align}

Our proof will make use of a criterion for strict positivity of the
density of a random variable, which was first established in the case
of finite dimensional diffusions by Ben Arous and L\'eandre
\cite{b:BenArousLeandre91b}.  It was then extended to general random
variables defined on Wiener space by Aida, Kusuoka and Stroock
\cite{b:AidaKusuokaStroock93}. We follow the presentation in Nualart
\cite{b:Nualart95a}. However, there is one major difference between
our case and the classical situation treated in those references. As
noted at the start of Section \ref{sec:RepMallMatrix}, our SPDE can be
solved pathwise. This means that the Wiener process
$W(t)=(\;W_k(t)\;)_{k\in \Zf}$ can be replaced by a fixed trajectory
in $\Ot\eqdef C([0,t]; \R^\Zf)$. Hence, we do not really need the
notion of a skeleton. Because of this we can prove a result which is
slightly more general than usual in that our controls need not belong
to the Cameron Martin space. Let $Q \in
\mathcal{L}(\R^\Zf;\L^2)$ be such that if $\{q_k , k\in\Zf\}$ is a
standard basis for $\R^\Zf$, then $Qq_k=\base_k$.

The main result of this section is the following:
\begin{theorem}
  \label{thm:positivity}
  Assume that $S_\infty=\L^2$.
  Let $t >0$, and $\Pi$ be an orthogonal projection of $\L^2$ onto a
  finite dimensional subspace $S \subset S_\infty=\L^2$. Let $x \in S$ be
  such that for some $0 < s < t$ and all $\sol \in \L^2$ there exists
  $H \in C([s,t];\R^\Zf)$ such that the solution of 
  \begin{equation}\label{eq:solH}
    \left\{
      \begin{aligned}
    &\frac{\partial \sol^H}{\partial r}(r) + B(\sol^H(r),\sol^H(r)) =
    \nu\Delta \sol^H(r) + Q \frac{\partial H}{\partial r}(r),\; r> s\\
    &\sol^H(s)=\sol
      \end{aligned}\right.
      \end{equation}
      satisfies $\Pi \sol^H(t)=x$. Then the density $p(t,\cdot)$
      of the law of the random variable $\Pi \sol(t)$ satisfies $p(t,x) >0$.
\end{theorem}
Note that equation \eqref{eq:solH} makes perfect sense even when $H$
is not differentiable in time.  To see this define
$\bar{\sol}^H(s)\eqdef\sol^H(s)-QH(s)$ and note that the equation for
$\bar{\sol}^H(s)$ is a well known type of equation, to which existence
and uniqueness results apply (cf. \cite{b:CoFo88,b:FoiasProdi67}).

\begin{corollary}\label{c:RussianControl} Let $\Pi$ be any
  orthogonal projection of $\L^2$ onto a finite dimensional space.  If
  $S_\infty = \L^2$ then the density $p(t,\;\cdot\;)$ of the random
  variable $\Pi \sol(t)$ satisfies $p(t,x)> 0$ for all $x \in \Pi
  \L^2$.
\end{corollary}
\begin{proof}[Proof of Corollary \ref{c:RussianControl}]
  The results almost follows from the published version of Theorem 9
  in Agrachev and Sarychev \cite{b:AcrachevSarychev03Pre} which states
  that under some assumptions for any $s \in (0,t)$, the
  controllability assumption of Theorem \ref{thm:positivity} is
  satisfied with a control $H \in W^{1,\infty}(s,t;\R^\Zf)$. The
  increasing family of sets which describes the way the randomness
  spreads is slightly different the $S_n$. Furthermore, they do not
  state the result for arbitrary projection only the span of a finite
  number of Fourier modes. However, in private communications with the
  authors they have verified that the sets $S_n$ may be used and that an
  arbitrary finite dimensional projection may be taken. As an aside, Romito
  \cite{b:Romito02Pre} has proven this formulation of controllability of the Galerkin
  approximations under our assumptions.
\end{proof}

\begin{remark}
  It is worth pointing out that the exact control ability of the
  projections if far from the exact controllability in all of
  $\L^2$. In fact the later does not hold with smooth in space and  $L^2$
  in time controls. This would imply that the density was supported on
  $\L^2$ which is not true as its support is contained in functions
  which are analytic is space \cite{b:Mattingly98b,b:Mattingly02a}.
\end{remark}

The rest of the section is devoted to the proof of Theorem
\ref{thm:positivity}. Our proof is based on the following result,
which is a variant of Proposition 4.2.2 in Nualart \cite{b:Nualart95a}.
\begin{proposition}\label{p:pos}
 Let
$F \in C( \Ot; S)$ where $S=\Pi \L^2$ is a finite dimensional vector
space such that $ H \mapsto F(H)$ is twice differentiable in the
directions of $H^1(0,t;\R^\Zf)$, and these exist $DF(\cdot) \in C\big(
\Ot; \big[L^2(0,t;S)\big]^{\Zf}\big)$ and $D^2F(\cdot) \in C\big(\Ot; \big[L^2(
(0,t)^2; S)\big]^{\Zf\times\Zf}\big)$ such that for all $j,\ell \in\Zf$ and
$h,g \in L^2(0,t;\R)$,
\begin{align*}
  \frac{d\ }{d \epsilon} F\Big( H + \epsilon\int_0^\cdot q_j h(s) ds
  \Big)\Big|_{\epsilon=0} & = \int_0^t D_{j,s}F(H)h(s)ds,\\
  \frac{d^2\ }{d\epsilon\;d\delta} F\Big(  H  + \epsilon\int_0^\cdot
  q_j h(s) ds+ \delta \int_0^\cdot
  q_\ell g(s) ds\Big)\Big|_{\substack{\epsilon=0\\\delta=0}} &= \int_0^t \int_0^t
  D^2_{j,s;\ell,r} F(H)h(s) g(r) ds\;dr.
\end{align*}
We assume moreover that 
\begin{equation}
  \label{eq:8.1}
  H \mapsto \big( F(H), DF(H), D^2F(H)\Big)
\end{equation}
is continuous and locally bounded from $\Ot$ into $S \times \big[
L^2(0,t; S)\big]^{\Zf} \times \big[L^2([0,t]^2;S)\big]^{\Zf\times\Zf }$, and that there exists $H
\in \Ot$  such that $F(H)=x$ and
\begin{equation}
  \label{eq:8.2}
  \det \Big( \sum_{k\in\Zf} \int_0^t D_{k,s} F^i(H) D_{k,s} F^j(H)
  ds\Big) >0\;.
\end{equation}
Then, if the law of $F(W)$ has a density $p(t,\cdot)$, $p(t,x)
>0$.
\end{proposition}
\begin{proof}
  Since the proof is almost identical to that of Proposition 4.2.2 in
  \cite{b:Nualart95a}, we only indicate the differences with
  the latter proof.

  In this proof, $H$ is one specific element of $\Ot$ satisfying
  $F(H)=x$ and \eqref{eq:8.2}. Let $h_j(s)=D_{\cdot,s} F^j(H)$,
  $j \in \Zf$. Clearly, $h_j \in L^2(0,t; \R^\Zf)$. For each $z \in
  \R^\Zf$, let 
  \begin{equation*}
    (T_zW)(t)\eqdef W(t) + \sum_{j \in \Zf} z_j \int_0^t h_j(s) ds \;,
  \end{equation*}
  and 
  \begin{equation*}
    g(z,W)=F(T_zW)-F(W)\;.
  \end{equation*}
It follows from our assumptions that for any $W \in \Ot$, $g(\cdot,W)
\in C^2(B_1(0);\R^\Zf)$, and for any $\beta >1$, there exists $C(\beta)$
such that
\begin{align*}
  \|W\|_{\infty,t} \leq C(\beta) \Rightarrow
  \|g(\cdot,W)\|_{C^2(B_1(0))} \leq \beta\;,
\end{align*}
where $\|W\|_{\infty,t}\eqdef \sup_{0\leq s\leq t} |W(s)|$, and the notation 
$B_\alpha(0)$
stands for the open ball in $\R^\Zf$ centered at 0, with radius $\alpha$.

Assume for a moment that in addition
\begin{align*}
  \big| \det g'(0) \big| \geq \frac1\beta \;.
\end{align*}
It then follows from Lemma 4.2.1 in \cite{b:Nualart95a} that
there exists $c_\beta \in (0,\frac1\beta)$ and $\delta_\beta >0$ such
that $g(\cdot,W)$ is  diffeomorphic from $B_{c_\beta}(0)$ onto a
neighborhood of $B_{\delta_\beta}(0)$.

We now define the random variable ${\mathcal H}_\beta$, which plays exactly the
same role in the rest of our proof as $H_\beta$ in \cite{b:Nualart95a}, 
but
is defined slightly differently. We let
\begin{align*}
 {\mathcal H}_\beta=k_\beta(\|W\|_{\infty,t}) \alpha_\beta(|\det \ip{DF^i(W)}{DF^j(H)}|)\;,
\end{align*}
where $\ip{\cdot}{\cdot}$ denotes the scalar product in
$L^2(0,t;\R^\Zf)$;
$ k_\beta,\alpha_\beta \in C ({\R};[0,1])$, $k_\beta(x)=0$
whenever $|x| \geq \beta$, $k_\beta(x)=1$ whenever 
$|x| \leq \beta -1$; $\alpha_\beta(x)=0$ whenever 
$|x|\leq 1/\beta$, $\alpha_\beta(x) > 0$ whenever 
$|x| > 1/\beta$, and $\alpha_\beta(x)=1$ if $|x| \ge 2/\beta$.

The rest of the proof follows exactly the argument in
\cite{b:Nualart95a}
 pages 181 and 182. We only have to make explicit the sequence 
$T_N^H,\:N=1,2, \ldots$ of absolutely continuous transformations of $\Ot$ 
equipped with Wiener measure, which is used at the end of the proof. 
For $N=1,2, \ldots$, let $t_i^N=(i\,t)/N$, $0 \leq i \leq N$. We define
\begin{equation*}
(T_N^H W)(t)=W(t)+ \int_0^t(\dot{H}_N(s)-\dot{W}_N(s))ds,
\end{equation*}
where
\begin{equation*}
\begin{aligned}
\dot{H}_N(s)&= \sum_{i=1}^{N-1} \frac{H(t_i^N)-H(t_{i-1}^N)}{t_i^N-t_{i-1}^N} 
{\bf 1}_{[t_i^N, t_{i+1}^N)}(s),\\
\dot{W}_N(s)&= \sum_{i=1}^{N-1} \frac{W(t_i^N)-W(t_{i-1}^N)}{t_i^N-t_{i-1}^N} 
{\bf 1}_{[t_i^N, t_{i+1}^N)} (s).
\end{aligned}
\end{equation*}
For any $W \in \Ot$, as $N \rightarrow \infty$,
\begin{equation*}
\sup_{0 \leq s \leq t} \left|\left(T_N^H W\right)(s) -H(s)\right| 
\rightarrow 0\;.
\end{equation*}
Hence, by the continuity of $F$ and $DF$, we also have
\begin{align*}
F(T_N^HW) &\rightarrow F(H)\\
DF(T_N^HW)&\rightarrow DF(H),
\end{align*}
and moreover
\begin{equation*}
\lim_{M\rightarrow\infty} \sup_N {\P}\left(\|T_N^H W\|_{\infty,t}>M\right)=0.
\end{equation*}
This provides exactly the version of (H2) from Nualart, which is needed here to complete the proof.
\end{proof}

All that remains is to prove the following lemma:
\begin{proposition}\label{pro8.5}
Under the assumptions of Theorem \ref{thm:positivity}, if $F=\Pi \sol(t)$, 
there exists $H \in \Ot$ such that $F(H)\eqdef \Pi \sol^H(t)=x$, 
and \eqref{eq:8.2} holds.
\end{proposition}
\bpf
Let $s \in (0,t)$ be the time which appears in the assumption of Theorem 
\ref{thm:positivity}. Since $S_\infty=\L^2$, it follows from \eqref{eq:noDeg} that
\begin{equation*}
{\P}\Biggl(\bigcap_{\phi\in \H^1, \; \phi\not=0}\left\{
\ip{{\mathcal M}(s)\phi}{\phi}>0\right\}\Biggr)=1.
\end{equation*}
We choose an arbitrary Brownian trajectory $W$ in this set of measure one. 
We then choose, according to the assumption of Theorem \ref{thm:positivity}, 
an $H_x \in \Ot$ satisfying $H_x(s)=0$, such that the solution 
$\{\sol^{H_x}(r),\: s \leq r \leq t\}$ of \eqref{eq:solH} 
with initial condition
$\sol^{H_x}(s)=\sol(s;\sol_0,W)$, the value at time $s$ of the solution of 
\eqref{eq:vorticity} corresponding to the trajectory $W$ which we chose above, 
satisfies $\Pi \sol^{H_x}(t)=x$. It remains to show that the condition 
\eqref{eq:8.2} 
is satisfied, with $H \in \Ot$ defined by
\begin{equation*}
H(r)=
\begin{cases} 
W(r), &\text{ if }0 \leq r \leq s;\\
W(s)+H_x(r), &\text{ if }s < r \leq t.
\end{cases}
\end{equation*}
We shall write
\begin{equation*}
\sol^H(r)=
\begin{cases}
\sol(r,W),&\text{ if }0\leq r \leq s,\\
\sol^{H_x}(r),&\text{ if }s < r \leq t.
\end{cases}
\end{equation*} 
All we need to show is that for any fixed $\phi \in S$, $\phi \not=0$,

\begin{equation*}
\sum_{k\in\Zf} \int_0^t\ip{D_{k,r}\sol^H(t)}{\phi}^2 dr >0.
\end{equation*}
Note that
\begin{equation*}
\begin{aligned}
\sum_{k \in \Zf}\int_0^t\ip{D_{k,r}\sol^H(t)}{\phi}^2 dr  
&\geq \sum_{k \in \Zf} \int_0^s\ip{D_{k,r}\sol^H(t)}{\phi}^2 dr\\
&=\sum_{k \in \Zf} \int_0^s \ip{D_{k,r} \sol^H(s)}{U^{t,\phi}(s)}^2 ds,
\end{aligned}
\end{equation*}
where $\{U^{t,\phi}(s),\:  0 \leq s \leq t\}$ solves the backward equation 
\eqref{2.2}, associated to the corresponding $\sol$ trajectory, and the 
last identity follows from Proposition \ref{prop2.3}. In view of our specific 
choice of $W$, since from Lemma \ref{l:JJmoments} $U^{t,\phi}(s) \in \H^1$, 
it suffices 
to check that $U^{t,\phi}(s)\not=0$. This follows from ``backward 
uniqueness'' (maybe we should say ``forward uniqueness'' since the 
$U^{t,\phi}$ equation is a backward equation !), see
 Theorem I.1 in Bardos, Tartar \cite{b:BardosTartar73}, whose assumptions are 
clearly verified by the $U^{t,\phi}(\cdot)$ equation
\eqref{2.2}. Specifically, since if $\phi=0$ then  $U^{t,\phi}(s)=0$
for all $s \leq t$ one knows that no other terminal condition $\phi$
at time $t$ can lead
to $U^{t,\phi}(s)=0$ for $s \leq t$.
\epf
\section{Conclusion}

We have proven under reasonable nondegeneracy conditions that the law
of any finite dimensional projection of the solution of the stochastic
Navier Stokes equation with additive noise possesses a smooth,
strictly positive density with respect to Lebesgue measure. In
particular, it was shown that four degrees of freedom are sufficient
to guarantee nondegeneracy.

It is reasonable to ask if four is the minimal size set which produces
finite dimensional projections with a smooth density. The
nondegeneracy condition concentrates on the wave numbers were both the
$\sin$ and $\cos$ are forced. Since this represents the translation
invariant scales in the forcing, it is a mild restriction to require
that whenever either of the $\sin$ or $\cos$ of a given wave number is
forced, then both are forced. Under this assumption, forcing only two
degrees of freedom corresponds to forcing both degrees of freedom
associated to a single wave number $k$. It is easy to see that the
subspace $\{ u \in \L^2 : \ip{u}{\sin(j\cdot x)}= \ip{u}{\cos(j\cdot
  x)}=0 \mbox{ for $j\neq k$}\}$ is invariant under the dynamics with
such a forcing. Hence if the initial  condition lies in this two dimensional
subspace, the conclusions of Theorem \ref{theo1.1} fail to hold. See
\cite{b:HairerMattingly04} for a more complete discussion of this and
other cases where the nondegeneracy condition fails.

We have concentrated on the 2D Navier Stokes equations forced by
finite number of Wiener processes. However there are a number of ways
one could extend these results. The choice of a forcing with finitely
many modes was made for simplicity in a number of technical lemmas, in
particular in section \ref{sec:ladderDef}.  There appears to be no
fundamental obstruction to extending the method to the cases with
infinitely many forcing terms if the covariances satisfy an
appropriate summability condition. In addition, the methods of this
paper should equally apply to other polynomial nonlinearities, such as
stochastic reaction diffusion equations with additive noise.  In
contrast, handling non-additive forcing in a nonlinear equation would
require nontrivial extensions of the present work. Since in the
linearization, stochastic integrals of nonadapted processes would
appear, it is not certain that the line of argument in this paper
would succeed.



\ack We would like to thank Andrei Agrachev, G\'erard Ben Arous, Yuri Bakhtin,
Weinan E, Thierry Gallouet, Martin Hairer, George Papanicolaou, and
Yakov Sinai for stimulating and useful discussions.  We also thank
Martin Hairer who after reading a preliminary version of this paper
suggested a better notation for the Fourier basis which we proceeded
to adopt. We also thank the anonymous referee for her (or his)
comments and corrections. JCM was partial supported by the Institute
for Advanced Study and the NSF under grants DMS-9971087 and
DMS-9729992.  JCM would also like to thank the Universit\'e de
Provence and the IUF which supported his stays in Marseilles during
the summers of 2002-2004 when this work was undertaken.

        

\appendices

\section{Estimates on the Enstrophy}
Define $\mathcal{E}_0=\sum_{k\in\Zf} \|\base_k\|^2$ where the $\base_k$ were the
functions used to define the forcing in \eqref{eq:W}.  In general we
define the $\alpha$-th spatial moment of the forcing to be
$\mathcal{E}_\alpha=\sum_{k\in\Zf} |k|^{2\alpha}\|e_k\|^2$. With this notation,
we have the following estimate.
\begin{lemma} \label{eq:EnstrophyControl}
Given any $\epsilon \in (0,1)$, there exists a
  $\gamma=\gamma(\epsilon)$ such that
  \begin{equation*}
  \P\left\{ \sup_{s \in [0,t]} \|\sol(s)\|^{2} +2\epsilon\nu
  \int_0^s \|\sol(r)\|_1^2 dr - \mathcal{E}_0 s >
  \|\sol(0)\|^2 +K \right\} < e^{-\gamma K}    
  \end{equation*}
  for all $K \geq 0$.
\end{lemma}
\bpf The Lemma follows from the exponential martingale estimate after
one notices that the quadratic variation of the martingale in the
equation for the enstrophy is controlled by
$\int_0^t\|\sol(r)\|_1^2\;dr$. See \cite{b:Mattingly02} Lemma A.2 or
\cite{b:Mattingly02} Lemma A.1 for the details and related lemmas in
exactly this setting or Lemma \ref{l:nablaSol} below for a similar
argument.\epf

{}From the previous result, we obtain the following.
\begin{corollary}\label{EnstrophyExpMoment} There exists a constant
  $\eta_0=\eta_0(T,\nu) >0$ so that for any $\eta \in (0,\eta_0]$
  there exists a constant $c=c(T,\eta,\nu)$ so that
  \begin{align*}
    \E\exp\Big( \eta \sup_{0\leq s \leq T} \| \sol(s)\|^2 \Big)
    \leq c \E \exp\left( \eta \|\sol(0)\|^2\right)
\intertext{and}
    \E\exp\Big( \nu \eta \int_0^T\| \sol(s)\|_1^2 ds \Big)
    \leq c \E \exp\left( \eta \|\sol(0)\|^2\right)
  \end{align*}
\end{corollary}

We will also need the following result which gives quantitative
estimates on the  regularization of the $\H^1$-norm. 
\begin{lemma} \label{l:nablaSol}Given a time $T > 0$ and $p\geq 0$,
  there exists positive constants $c=c(\nu,T,p,\mathcal{E}_1)$ so that  
  \begin{align*}
    \E \sup_{0\leq s \leq T} \big[\|\sol(s)\|^{2} +s \nu
    \|\sol(s)\|_1^{2}\big]^p \leq c[ 1 + \|\sol(0)\|^{4p}]
  \end{align*}
\end{lemma}
\bpf
Defining $\zeta(s)= \|\sol(s)\|^2+ s \nu \|\sol(s)\|_1^2$, we have for all $s \in
[0,T]$
\begin{multline*}
  \zeta_s  -\mathcal{E}_0s - \frac12 \mathcal{E}_1s^2 = \|\sol(0)\|^2
  + \int_0^s 2r \nu \big[ -\nu \|  \sol(r)\|_2^2 + \ip{B(\sol(r),\sol(r))}{\Delta \sol(r)}\big]dr  \\
  - \nu \int_0^s \|\sol(r)\|_1^2 dr + \sum_{k\in\Zf} \int_0^s 2(1+
  r\nu |k|^2)\ip{\sol(r)}{\base_k} dW_k(r) \ .
\end{multline*}
Since
\begin{align*}
  |\ip{B(\sol,\sol)}{\Delta \sol}| \leq c\|\sol\|_\frac12
  \|\sol\|_1 \|\sol\|_2 \leq \nu \|\sol\|_2^2 +
  c(\nu)\|\sol\|^4,
\end{align*}
one has
\begin{equation}\label{eq:solRegulation}
  \sup_{0\leq s \leq T} \zeta_s \leq c(1+T^2) +
  cT^2\sup_{0\leq s \leq T} \|\sol(s)\|^4 + \sup_{0\leq s \leq T} N_s
\end{equation}
where 
\begin{align*}
N_s&= -\nu \int_0^s [\|\sol(r)\|^2+ r\|\sol(r)\|_2^2] dr + M_s 
\\ \intertext{and}
M_s&=
\sum_{k \in \Zf} \int_0^s 2(1+ r\nu |k|^2)\ip{\sol(r)}{\base_k} dW_k(r)\;.  
\end{align*}
Notice that for all $s \in [0,T]$ and $\alpha>0$ sufficiently small
$N_s \leq M_s - \frac\alpha2 [M,M]_s$, where $ [M,M]_s$ is the
quadratic variation of the martingale $M_s$. Hence the exponential
martingale estimate implies $\P( \sup_{s\leq T} N_s > \beta ) \leq \P(
\sup_{s\leq T} M_s - \frac\alpha2 [M,M]_s> \beta ) \leq \exp(-\alpha
\beta)$.  This implies that the last term in \eqref{eq:solRegulation}
can be bounded by a constant depending only on $\alpha, T,\nu$ and the
power $p$. By Corollary \ref{EnstrophyExpMoment}, the third term in
\eqref{eq:solRegulation} can be bounded by a constant which depends on
the initial condition as stated as well as $\alpha, T,\nu$ and the
power $p$.
\epf

\section{Estimates on the Linearization and its Adjoint}
\label{sec:VUmoments}

Define the action of linearized operator $J_{s,t}$ on a $\phi \in
\L^2$ by
\begin{equation}
  \label{eq:J}
\left\{
\begin{aligned}
  &\frac{\partial \ }{\partial t} J_{s,t}\phi = \nu \Delta  J_{s,t}\phi
  + B(\sol(t),  J_{s,t}\phi) +  B( J_{s,t}\phi,\sol(t)) \qquad   0
  \leq s \leq t;\\
   &J_{s,s}\phi= \phi
\end{aligned}
\right.
\end{equation}
and its time reversed, $\L^2$-adjoint $\bJ_{s,t}^*$  acting on $\phi \in
\L^2$ by
\begin{equation}
  \label{eq:J*}
\left\{
\begin{aligned}
  &\frac{\partial \ }{\partial s} \bJ_{s,t}^*\phi + \nu \Delta  \bJ_{s,t}^*\phi
  + B(\sol(s),  \bJ_{s,t}^*\phi) - C( \bJ_{s,t}^*\phi,\sol(s))=0 \qquad   0
  \leq s \leq t;\\
  &\bJ_{t,t}^*\phi= \phi
\end{aligned}
\right.
\end{equation}

If we define the operator $Q:\R^{\Zf} \rightarrow \L^2$ by
$(x_k)_{ k \in\Zf} \mapsto \sum x_k \base_k$, then $V_{k,s}(t)=J_{s,t}Qq_k$
were $\{q_k:  k \in\Zf\}$ is the standard basis for $\R^{\Zf}$ and
$U^{t,\phi}(s)=\bJ^*_{s,t}\phi$. Similarly for $h \in L^2_{loc}(\R^{\Zf}_+$), 
$D_h\sol(t)=\int_0^t J_{s,t}Qh(s) ds$.
\begin{lemma}
  \label{l:JJmoments} For any $T_0 > 0$,  $\eta >0$ and $\alpha \in\{0,1\}$ there exists
  constants $\gamma=\gamma(\nu,\eta,T_0^\alpha)$ and
  $c=c(\nu,T_0^\alpha,\eta)$ so that for all $\phi \in \L^2$
  and $T \leq T_0$
  \begin{align*}
    \sup_{0\leq s \leq t \leq T} \|J_{s,t}\phi\|^2+ \gamma
    (t-s)^\alpha\|J_{s,t}\phi\|_1^2 &\leq \exp\Big(\eta \int_0^T
    \|\sol(r)\|_1^2 dr  + cT  \Big)\Big(\|\phi\|^2
    +(1-\alpha)\|\phi\|_1^{2}\Big)\\ 
    \sup_{0\leq s \leq t \leq T} \|\bJ_{s,t}^*\phi\|^2+ \gamma
    (t-s)^\alpha\|\bJ_{s,t}^*\phi\|_1^2 &\leq \exp\Big(\eta \int_0^T
    \|\sol(r)\|_1^2 dr + cT\Big)\Big(\|\phi\|^2 +
    (1-\alpha)\|\phi\|_1^{2}\Big)\\ 
  \end{align*}
where on the right hand side $  (1-\alpha)\|\phi\|_1^{2}=0$ when
$\alpha=1$  by convention for all $\phi$ even if $\|\phi\|_1=\infty$.
\end{lemma}
\bpf We start by deriving a number of bounds on the nonlinear
terms. Using Lemma \ref{l:Bestimates} and standard interpolation
inequalities produces for any $\delta >0$ and $\eta >0$ and some $c$,
\begin{align*}
  2|\ip{B(\phi,\sol)}{\phi}| &\leq c \|\phi\|^\frac32
  \|\phi\|_1^\frac12 \|\sol\|_1 \leq \delta \|\phi\|_1^2 + \|\phi\|^2
  \Big( \eta \|\sol\|_1^2 + \frac{c}{\delta \eta^2}
  \Big)\\
  2|\ip{B(\phi,\sol)}{\Delta \phi}| &\leq c\|\sol\|_1
  \|\phi\|_1^\frac12 \|\phi\|^\frac12 \|\phi\|_2 \leq \delta
  \|\phi\|_2^2 + \eta \|\sol\|_1^2 \|\phi \|_1^2 + \frac{c}{\eta
    \delta^2} \|\sol\|_1^2 \|\phi \|^2\\
  2|\ip{B(\sol,\phi)}{\Delta \phi}| &\leq c\|\sol\|_1 \|\phi\|_1
  \|\phi\|^\frac14 \|\phi\|_2^\frac34 \leq \frac{c}{(\eta
    \delta)^\frac32} \|\sol\|_1 \|\phi\|_1 \|\phi\| + (\eta
  \delta)^\frac12\|\sol\|_1 \|\phi\|_1\|\phi\|_2\\
  &\leq \delta \|\phi\|_1^2 + \frac{c}{\eta^3\delta^4}\|\sol\|_1^2
  \|\phi\|^2 + \frac{\delta}2\|\phi\|_2^2 + \eta \|\sol\|_1^2
  \|\phi\|_1^2\\
  2|\ip{B(\Delta \phi,\sol)}{\phi}| &\leq c\| \phi\|_{1+\epsilon}
  \|\phi\|_{1-\epsilon} \|\sol\|_1 \leq c \|\sol\|_1
  \|\phi\| \|\phi\|_2 \leq \delta \|\phi\|_2^2 +
  \frac{c}{\delta}\|\sol\|_1^2 \|\phi\|^2
\end{align*}
We begin with bounding $J$ expression with $\alpha=0$. Setting $\zeta_r=\|J_{s,s+r}\phi\|^2 +
\gamma \|J_{s,s+r}\phi\|_1^2$, $\sol_r=\sol(s+r)$ and
$J_r=J_{s,s+r}\phi$  and using the above estimates with appropriately
chosen constants produces
\begin{align*}
  \frac{\partial \zeta_r}{\partial r} = & \frac{\partial\ 
  }{\partial r} \|J_{r}\|^2 + \gamma \frac{\partial \ }{\partial r} \|J_{r}\|_1^2\\
  =& - 2\nu \|J_{r}\phi\|_1^2 + 2\ip{B(J_r,\sol)}{J_r} \\ &
  \qquad \qquad + 2\gamma\big[- \nu \|J_{r}\|_2^2
  +\ip{B(J_r,\sol)}{\Delta J_r} +
  \ip{B(\sol,J_r)}{\Delta J_r} \big]\\
  \leq & \Big( \frac{c}{\eta^2\nu} + \frac{\eta}{2} \|\sol_r\|_1^2 +
  \gamma c\big(\frac{1}{\eta\nu^2} + \frac{1}{\eta^3\nu^4} \big)\|\sol_r\|_1^2 \Big)\|J_r\|^2 + \gamma\eta\|\sol_r\|_1^2 \| J_r\|_1^2
\end{align*}
Hence for $\gamma$ sufficiently small there exist a constant
$c$ so for all $r > 0$
\begin{align*}
  \frac{\partial \zeta_r}{\partial r} \leq [ \frac{c}{\eta^2\nu} + \eta
  \|\sol_r\|_1^2 ]\zeta_r
\end{align*}
which proves the first result for $J_{s,t}$. The result for $\alpha=1$
is identical except that we take $\zeta_r=\|J_{s,s+r}\phi\|^2 +
\gamma r \|J_{s,s+r}\phi\|_1^2$ and hence there is an extra term from the
differentiating the coefficient of $\|J_{s,s+r}\phi\|_1^2$ and
the fact that the resulting constants depend on the time interval
$[0,T]$ over which $r$ ranges. See the proof of Lemma \ref{l:nablaSol}.

Turning to $\bJ_{s,t}^*$, we set $\bar \zeta_r=\|\bJ_{t-r,t}^*\phi\|^2 +
\gamma \|\bJ_{t-r,t}^*\phi\|^2$, $\bar \sol_r=\sol(t-r)$ and $\bar
J_r^*=\bJ_{t-r,t}^*\phi$ for all $r \in [0,t]$. The bar is to remind us
that the process is time reversed. The argument proceeds as in the
previous case. Again using the estimates above, we obtain
\begin{align*}
    \frac{\partial \bar \zeta_r}{\partial r} \leq \Big( \frac{c}{\nu(1+ \eta^2)}
    + \big( \frac{\eta}{2} + \frac{\gamma c}{\nu} \big) \|\sol_r\|_1^2 \Big) \|\bJ_r^*\|^2 + \gamma \eta \|\sol_r\|_1^2
    \|\bJ^*_r\|_1^2 \ .
\end{align*}
Hence for $\gamma$ small enough and $r \in [0,t]$,
\begin{align*}
  \frac{\partial \bar \zeta_r}{\partial r} \leq \Big[ \frac{c}{\nu(1+
    \eta^2)} + \eta \|\sol_r\|_1^2 \Big]\bar \zeta_r
\end{align*}
which proves the result.
\epf

\section{Higher Malliavin Derivatives of $\sol(t)$}
\label{sec:HigherMalliavainD}
For notational brevity define $\tilde B(f,g)= B(f,g)+B(g,f)$ for $f,g
\in \L^2$. For $k\in\Zf$, we  define $J_{s,t}^{(1)}e_k=
J_{s,t}e_k$. For $n>1$ define $J_{s,t}^{(n)}$ acting on
$\phi=(e_{k_1},\dots,e_{k_n})$ with
$k_1,\dots,k_n\in\Zf$ and with time
parameters $s=(s_1,\dots,s_n) \in \R_+^n$ by the equations
\begin{align*}
  \frac{\partial \ }{\partial t}J_{s;t}^{(n)}\phi &= \nu \Delta J_{s;t}^{(n)}\phi
  + \tilde B(\sol(t),  J_{s;t}^{(n)}\phi) + F_{s;t}^{(n)}\phi; \quad  t>
  \vee s   \\
  J_{s;t}^{(n)} & = 0; \quad t \leq \vee s
\end{align*}
where $\vee s=s_1 \vee \cdots \vee s_n$. The operator $F_{s,t}^{(n)}$
applied to $\phi$ is defined by
\begin{align*}
  F_{s;t}^{(n)}\phi =\sum_{(\alpha,\beta)\in \mbox{\tiny part}(n)} \tilde
  B\big(J_{s_\alpha;t}^{(|\alpha|)}\phi_\alpha, 
    J_{s_\beta;t}^{(|\beta|)}\phi_\beta\big)\ .
\end{align*}
Here $\mbox{part}(n)$ is the set of partitions of $\{1,\dots,n\}$ into
two sets, neither of them empty. $|\alpha|$ is the number of elements
in $\alpha$,
$\phi_\alpha=(\phi_{\alpha_1},\dots,\phi_{\alpha_{|\alpha|}})$ and
$s_\alpha=(s_{\alpha_1},\dots,s_{\alpha_{|\alpha|}})$. The partition
$(\alpha,\beta)$ and $(\beta,\alpha)$ are viewed as the same
partition.

First observe that when $n=1$, Lemma \ref{l:JJmoments} says that for any $\eta >0$ there is
a $c=c(T,\eta)$ so for all $\phi \in \R^{\Zf}$
\begin{align*}
\sup_{0\leq s \leq t\leq T}\|J_{s;t}^{(1)}e_k\|_1&\leq  c
 \exp\Bigl( \eta\int_0^T \|\sol(r)\|_1^2 dr
 \Bigr) \leq  c \exp\Bigl( \eta\int_0^T \|\sol(r)\|_1^2 dr
 \Bigr)\ .  
\end{align*}

For $n > 1$ with again $\phi=(e_{k_1},\dots,e_{k_n})$ and $s \in
\R_+^n$, we have the following estimate on $F^{(n)}$
\begin{align*}
  \|F_{s;t}^{(n)}\phi \|\leq & \sum_{(\alpha,\beta)\in \mbox{\tiny
      part}(n)} \| \tilde
  B\big(J_{s_\alpha;t}^{(|\alpha|)}\phi_\alpha,
  J_{s_\beta;t}^{(|\beta|)}\phi_\beta\big)\| \leq c
  \sum_{(\alpha,\beta)\in \mbox{\tiny part}(n)} \|
  J_{s_\alpha;t}^{(|\alpha|)}\phi_\alpha\|_1\|
  J_{s_\beta;t}^{(|\beta|)}\phi_\beta\|_1\ .
\end{align*}
Then the variation of constants formula implies that $J_{s;t}^{(n)}\phi =
\int^t_{\vee s} J_{r,t} F_{s;r}^{(n)} \phi dr $ and hence
\begin{align*}
  \|J_{s;t}^{(n)}\phi\|_1 &\leq \int^t_{\vee s} \|J_{r,t}
  F_{s;r}^{(n)} \phi\|_1 dr\\
  & \leq c\exp\Bigl( \eta\int_0^T \|\sol(r)\|_1^2 dr \Bigr)
  \Bigl(\int^t_{\vee s} \frac{1}{(r-\vee
    s)^\frac12}dr\Bigr)\sup_{\substack{\tau=(\tau_1,\cdots,\tau_n)\\
      0\leq \tau_i\leq t\\ \vee \tau \leq r \leq
      t}}\|F_{\tau;r}^{(n)}\phi\|\\
  &\leq c\exp\Bigl( \eta\int_0^T \|\sol(r)\|_1^2 dr
  \Bigr)\sum_{(\alpha,\beta)\in \mbox{\tiny part}(n)} \sup_{\tau,r}\|
  J_{\tau_\alpha;r}^{(|\alpha|)}\phi_\alpha\|_1\sup_{\tau,r}\|
  J_{\tau_\beta;r}^{(|\beta|)}\phi_\beta\|_1
\end{align*}
Proceeding inductively, one obtains for any $\eta >0$ the existence
for a $c=c(T,\eta,n)$ and $\gamma=\gamma(n)$ so that  
\begin{align*}
  \|J_{s;t}^{(n)} \phi\|_1 \leq \int_s^t \| J_{r,t}
  F_{s,r} \phi\|_1 dr \leq c\exp\Bigl( \gamma\eta\int_0^T \|\sol(r)\|_1^2 dr
 \Bigr) \ .
\end{align*}
Now with $\phi$ and $s$ as above, 
\[
D^n_{s_1,k_1;\dots;s_n,k_n}\sol(t)=J^{(n)}_{s,t}\phi,
\]
and
since $\eta$ was a arbitrary positive constant, by redefining it, one
obtains that for any $\eta>0$, there
exits a $c=c(T,\eta,n)$ so that 
\begin{equation*}
  \sup_{ t\in (0,T]} \|D_{s_1,k_1;\dots;s_n,k_n}^{n}\sol(t)\|_1
  \leq  c\exp\Bigl( \eta\int_0^T \|\sol(r)\|_1^2 dr \Bigr).
 \end{equation*}
Combining this estimate with Lemma \ref{EnstrophyExpMoment}, we obtain
the following result. Letting
$\DD^\infty(\H_1)$ is the space of random variable taking values in
$\H_1$ which are inifintely differentiable in the Malliavin sense and
such that those derivatives have all moments finite, we have that: (see
\cite{b:Nualart95} page 62 for the definition of $\DD^\infty$)
\begin{lemma}\label{l:Dinfty}
  For any $\eta >0$, $t>0$, $p \geq 1$ and $n \geq 1$ there exists a constant
  $c=c(t,\nu,\eta,p,n)$ so that 
  \begin{align*}
    \E \left[\left(\sum_{k_1,\dots,k_n\in\Zf}\int_0^t\cdots\int_0^t \|
        D_{s_1,k_1;\dots;s_n,k_n}^{n}\sol(t)\|_1^2ds_1\cdots ds_n
       \right)^{p/2}\right] 
      \leq c \exp\left(\eta\|\sol(0)\|^2\right)\;.
    \end{align*}
Hence $\sol(t) \in \DD^\infty(\H_1)$ for all $t>0$. 
\end{lemma}

\section{Estimates on the Nonlinearity}

In the following, Lemma we collect a few standard estimates on the
nonlinearity and derive a few consequences from them.
\begin{lemma}  \label{l:Bestimates} Let $\alpha_i \geq 0$ and either
  $\alpha_1+\alpha_2+\alpha_3 > 1$ or both $\alpha_1+\alpha_2+\alpha_3=1$
  and $\alpha_i \neq 1$ for any $i$. Then the following estimates hold
  of all $f,g,h \in \L^2$ if
  the right hand side is well defined.
  \begin{equation*}
    |\ip{B(f,g)}{h}| \leq c \|f\|_{\alpha_1-1}\|g\|_{\alpha_2+1}\|h\|_{\alpha_3}
  \end{equation*}
 In addition we have the following estimate. For any $\epsilon >0$ there
 exists a $c=c(\epsilon)$ with 
  \begin{equation*}
    |\ip{\nabla B(f,g)}{\nabla g}| \leq c
     \|f\|_1\|g\|_1\|g\|_{1+\epsilon}
  \end{equation*}
\end{lemma}
\bpf

For the first result see Proposition 6.1 of \cite{b:CoFo88} and recall that
our $B$ is slightly different than theirs and that $\|\mathcal{K} f\|_1=
\|f\|_0$. After translation, the result follows. For the second result
we need to rearrange things. Setting $u=(u_1,u_2)= \mathcal{K} f$ we have
\begin{align*}
  |\ip{\nabla B(f,g)}{\nabla g}|&= |\int_{\T^2} \nabla [ ( u \cdot
   \nabla)g] \cdot \nabla g dx | 
\end{align*}
Observe that $ \nabla [ ( u \cdot \nabla)g] \cdot \nabla g= [(u \cdot
\nabla)\nabla g]\cdot \nabla g + (\nabla u \nabla g)\cdot \nabla g $.
Because $\nabla \cdot u=0$, $ [(u \cdot \nabla)\nabla g]\cdot \nabla
g= \frac12 \sum_i \nabla \cdot (u (\frac{\partial g}{\partial x_i})^2
)$. Hence, the integral of the first term is zero by stokes theorem
and the fact that we are on the torus.
   
The integral of the second term over the domain is made of a finite
number of terms of the form $\int \frac{\partial u_j}{\partial
  x_i}\frac{\partial g}{\partial x_j}\frac{\partial g}{\partial x_i}
dx $. This term is dominated by $| \frac{\partial u_j}{\partial
  x_i}|_{L^r}|\frac{\partial g}{\partial x_j}|_{L^p} |\frac{\partial
  g}{\partial x_i} |_{L^q}$ for any $r,p,q >1$ with
$\frac1r+\frac1q+\frac1p=1$. Recall that in two dimensions, the
Sobolev space $W^{1,2}$ is
embedded in  $L^r$ for any $r < \infty$. Hence by
taking $p=\frac12$, $q>2$ sufficiently close to 2 and $r$
correspondingly large, we obtain the bound $c\| \frac{\partial
  u}{\partial x_i}\|_1 \|\frac{\partial g}{\partial x_j}\|
\|\frac{\partial g}{\partial x_i} \|_\epsilon$ for any $\epsilon >0$
and some $c=c(\epsilon)$. The estimate is in turn bounded by $c\|f\|_1
\|g\|_1 \|g\|_{1+\epsilon}$.  \epf

\section{Lipschitz and Supremum Estimates}
\label{sec:LipschitzSupremumEstimates}
Let $S$ be a subspace of $\L^2$ spanned by a finite number of
$\cos(x\cdot k)$ and $\sin(x\cdot k)$. Let $\Pi$ be the orthogonal
projection onto $S$. Also let $\Pi_0$ be the orthogonal projection
onto the directions directly forced by Wiener processes as defined in
Proposition \ref{prop:formXYW}. 

Recall from Section \ref{sec:noSmallEVs} the definitions of
$\hol_{\alpha,[a,b]}(f)$, $\|f\|_\infty$, $\|f\|_{\infty,[a,b]}$
applied to functions of time taking values in $\L^2$.

For $0\leq s < t \leq T$ one has
\begin{align*}
  \Pi_0^\perp \sol(t) = \Pi_0^\perp e^{\nu \Delta(t-s)} \sol(s) -
  \Pi_0^\perp \int_s^t e^{\nu \Delta(t-r)} B(\sol(r),\sol(r)) dr .
\end{align*}
Since $\int f(x)dx=\int g(x)dx=0$, $\| B(f,g)\| \leq c \|\nabla f\|
\|\nabla g\|$ and $\|(e^{\nu \Delta (t-s)}-1)f\|\leq (1-e^{-\nu \lambda
  (t-s)})\|f\|$ for some fixed $\lambda>0$, one has
\begin{align*}
  \|\Pi_0^\perp[\sol(t)-\sol(s)] \| &\leq 2\Big[1-e^{-\lambda
    \nu(t-s)}\Big]\|\sol(s)\| + 2 \sup_{0\leq r\leq T} \|\nabla
  \sol(r)\|^2|t-s| \\&\leq c|t-s| \Big[ 1+ \|\nabla \sol\|_\infty^2
  \Big] \ .
\end{align*}
Similarly,
\begin{align*}
  \|U^{T,\phi}(t)-U^{T,\phi}(s)\| &\leq c|t-s| \|U^{T,\phi}(t)\| +
  c\|\nabla \sol\|_\infty \|\nabla U^{T,\phi}\|_\infty|t-s| \\
  & \leq c\big[ 1+ \|\nabla \sol\|_\infty^2 + \|\nabla
    U^{T,\phi}\|_\infty^2 \big]|t-s| \ .
\end{align*}
Also if 
\begin{align*}
  R(t)=\Pi_0 w(0) + \int_0^t \nu\Delta \Pi_0 \sol(r) - \Pi_0
    B(\sol(r),\sol(r)) dr
\end{align*}
then 
\begin{align*}
  \|R(t)-R(s)\| \leq c \int_s^t  \|\sol(r)\| dr + c\int_s^t
  \|\sol(r)\|^2 dr \leq c[ 1+ \|\sol\|_\infty^2] |t-s|
\end{align*}
Next observe that $\| \Pi B(f,g) \| \leq c\|f\| \|g\|$ and $\|\Pi
C(f,g)\| \leq c[ \|f\| \|\nabla g\| \wedge  \|g\|\|\nabla f\|]$. Thus
\begin{align*}
  \|\Pi C(f(t),g(t))- \Pi C(f(s),g(s)) \|& \leq \|\Pi
  C(f(t)-f(s),g(t))\| + \|\Pi C(f(s),g(t)-g(s)) \| \\&\leq c[
  \|\nabla g\|_\infty \|f(t)-f(s)\| + \|\nabla f\|_\infty
  \|g(t)-g(s)\| ]
\end{align*}
and 
\begin{align*}
   \|\Pi B(f(t),g(t))- \Pi B(f(s),g(s)) \|& \leq c[ \N{g}_\infty
   \|f(t)-f(s)\| +  \|f\|_\infty
   \|g(t)-g(s)\| ]
\end{align*}

We combine these observations in the following Lemma.
\begin{lemma}\label{l:auxLip} Let $\Pi$ and $\Pi_0$ be as above. In
  the notation of 
  \eqref{eq:DefLip} and \eqref{eq:defInf},
  \begin{align*}
    \hol_1(\Pi_0^\perp \sol) &\leq c\Big[ 1+ \|\nabla \sol\|_\infty^2
    \Big]\\
    \hol_1(U^{T,\phi}) & \leq c\big[ 1+ \|\nabla \sol\|_\infty^2 +
    \|\nabla U^{T,\phi}\|_\infty^2  \big]\\
    \hol_1(R) & \leq c[ 1+ \|\sol\|_\infty^2] \\
    \hol_1(\Pi B(f,g)) &\leq c [ \|f\|_\infty \hol_1(g) + \|g\|_\infty
    \hol_1(f)]\\
    \hol_1(\Pi C(f,g)) &\leq c [ \|\nabla f\|_\infty \hol_1(g) +
    \|\nabla g\|_\infty \hol_1(f)]
  \end{align*}
for all $f,g \in \L^2$ smooth enough that each term on the right
handside is finite.
\end{lemma}

Lastly we specialize these estimates to the setting of Proposition
\ref{prop:formXYW}. Let $X^\phi$, $R$, and $Y_k^{\phi}$ be as defined
in Proposition \ref{prop:formXYW} for $\sol(t)$ and $U^{T,\phi}$on the
interval $[0,T]$. Define $\chi^\phi= \Pi X^\phi$ and
$\Upsilon_k^{\phi}=\Pi Y_k^{\phi}$. We wish to obtain control of the
Lipschitz constants over an interval $[t,T]$ with $t \in (0,T)$.

Using the estimates from Lemma \ref{l:auxLip} and the fact that
$\|\nabla R\|_{\infty,[t,T]} \leq c\| R\|_{\infty,[t,T]} \leq c(
1+ \|\sol\|_{\infty,[t,T]}^2)$, we obtain
\begin{align*}
  \hol_{1,[t,T]}(\chi^\phi) \leq & c \big[ \hol_{1,[t,T]}(U^{T,\phi}) +
  \|\nabla \sol\|_{\infty,[t,T]} \hol_{1,[t,T]}(U^{T,\phi}) + \|\nabla
  U^{T,\phi}\|_{\infty,[t,T]} \hol_{1,[t,T]}(\Pi_0^\perp \sol) \\ & \qquad +
  \|\nabla R\|_{\infty,[t,T]} \hol_{1,[t,T]}(U^{T,\phi}) + \|\nabla
    U^{T,\phi}\|_{\infty,[t,T]}
  \hol_{1,[t,T]}(R)\big]\\
  \leq & c[1 + \|\nabla \sol\|_{\infty,[t,T]}^4 + \|\nabla
    U^{T,\phi}\|_{\infty,[t,T]}^4]
\end{align*}
and 
\begin{align*}
  \hol_{1,[t,T]}(\Upsilon_k^{\phi}) &\leq c \big[
  \hol_{1,[t,T]}(U^{T,\phi}) +\|
    U^{T,\phi}\|_{\infty,[t,T]}^2\big]\\&\leq c \big[ 1+ \|\nabla
    \sol\|_{\infty,[t,T]}^2 + \|\nabla U^{T,\phi}\|_{\infty,[t,T]}^2
  \big]\;.
\end{align*}
Similarly we have
\begin{align*}
  \|\chi^\phi\|_{\infty,[t,T]} &\leq c [ \|U^{T,\phi}\|_{\infty,[t,T]} +
  \|\nabla U^{T,\phi}\|_{\infty,[t,T]} \|\nabla \sol\|_{\infty,[t,T]}
  +\|\nabla U^{T,\phi}\|_{\infty,[t,T]}
  \|R\|_{\infty,[t,T]} ] \\
  & \leq c [ 1+ \|\nabla U^{T,\phi}\|_{\infty,[t,T]}^2 + \|\nabla
    \sol\|_{\infty,[t,T]}^2]
\end{align*}
and $\|\Upsilon_k^{\phi}\|_{\infty,[t,T]}\leq c \|U^{T,\phi}\|_{\infty,[t,T]}
$. 

Recall that $\N{f}_{\alpha,[t,T]}=\max( \|f\|_{\infty,[t,T]},
\hol_{\alpha,[t,T]}(f))$.  Combining the above estimates with Lemma
\ref{l:JJmoments} produces
\begin{align*}
  \N{\chi^\phi}_{1,[t,T]},
  \N{\Upsilon_k^{\phi}}_{1,[t,T]}  &\leq c[ 1+ \|\nabla U^{T,\phi}\|_{\infty,[t,T]}^4 
  + \|\nabla \sol\|_{\infty,[t,T]}^4]\\
  &\leq c\big[ 1+ \|\nabla\sol\|_{\infty,[t,T]}^{4}+ \exp\big(\eta
  \int_0^T \|\sol(r)\|^2_1 dr\big)\big]
\end{align*}
for all indices $i$, and $\phi$ with $\|\nabla \phi
\| \leq M$. Here $\eta >0$ is arbitrary but $c$ depends on the choice
of $\eta$ and $M$. In light of Corollary \ref{EnstrophyExpMoment} and
Lemma \ref{l:nablaSol} which control the right hand side, we obtain
the following result.
\begin{lemma}\label{l:XYMomentsForSNS} Given an $M>0$ define $S(M)=\{
  \phi : \|\nabla \phi \| \leq M \}$. Then for any $T >
  0$, $t \in (0,T)$, $p \geq 1$, and $\eta >0$ there exists a positive
  constant $c=c(\eta,p,t,\nu,\mathcal{E}_1,M,T)$ such that
  \begin{align*}
    \E \Biggl( \sup_{\phi \in S(M)} \Bigl[\N{ \chi^\phi}_{1,[t,T]}^p +
    \sup_{k \in \Zf} \N{\Upsilon_k^{\phi}}_{1,[t,T]}^p\Bigr] \Biggr) \leq c \exp
    \Bigl( \eta \|\sol(0)\|^2 \Bigr)
   \end{align*}
\end{lemma}


\frenchspacing
\bibliographystyle{plain}



                                

\def\cprime{$'$}

\end{document}